\documentclass[twoside]{article}
\usepackage{eurosym}
\usepackage[utf8]{inputenc}
\usepackage{graphicx}
\usepackage{latexsym}
\usepackage{amsmath}
\usepackage{amsfonts}
\usepackage{amssymb}
\usepackage{hyperref}
\usepackage[usenames]{color}

\newcommand{\displ}{\displaystyle}
\newcommand{\cA}{\mathcal{A}}
\newcommand{\cC}{\mathcal{C}}
\newcommand{\rC}{{\rm C}}
\newcommand{\rCh}{\widehat{\rm C}}
\newcommand{\rD}{{\rm D}}
\newcommand{\rE}{{\rm E}}
\newcommand{\cF}{\mathcal{F}}
\newcommand{\bff}{{\bf f}}
\newcommand{\bF}{{\bf F}}
\newcommand{\rG}{{\rm G}}
\newcommand{\rH}{{\rm H}}
\newcommand{\rI}{{\rm I}}
\newcommand{\J}{\mathbb{J}}
\newcommand{\rK}{{\rm K}}

\newcommand{\rL}{{\rm L}}
\newcommand{\rM}{{\rm M}}
\newcommand{\cP}{\mathcal{P}}
\newcommand{\rP}{{\rm P}}
\newcommand{\bP}{{\bf P}}
\newcommand{\rQ}{{\rm Q}}
\newcommand{\R}{{\mathbb{R}}}
\newcommand{\rR}{{\rm R}}
\newcommand{\rS}{{\rm S}}
\newcommand{\rT}{{\rm T}}
\newcommand{\cU}{\mathcal{U}}
\newcommand{\rV}{{\rm V}}
\newcommand{\rVh}{\widehat{\rm V}}
\newcommand{\bV}{{\bf V}}
\newcommand{\rW}{{\rm W}}
\newcommand{\xh}{\widehat{x}}
\newcommand{\by}{{\bf y}}

\newcommand{\fin}{\hfill\mbox{$\quad{}_{\Box}$}}
\newcommand{\fineq}{\vspace{.1cm}$\fin$}

\def\derp#1{\displ{\overset{\bullet}{#1}}}
\def\n#1{\displ{\|#1\|}}
\def\pe#1#2{\displ{\langle #1,#2\rangle }}

\newcommand{\proof}{{\sc Proof}.~}

\newtheorem{theo}{\bf \sffamily Theorem}
\newtheorem{prop}{\bf \sffamily Proposition}
\newtheorem{rem}{\bf \sffamily Remark}
\newtheorem{lemma}{\bf \sffamily Lemma}

\begin{document}
\title{\sffamily \bfseries A bang-bang optimal control for a nonlinear system modeling the Gate Control Theory of Pain}
\author{G. D\'{\i}az  and J.I. Díaz}  
\date{}

\maketitle
\begin{abstract}\small 
We consider a nonlinear system of coupled ordinary differential equations (representing the excitatory, inhibitory, and T-cell potentials) based on the Gate Control Theory of Pain, initially proposed by R. Melzack and P.D. Wall in 1965, and later mathematically modeled by N.F. Britton and S.M. Skevington in 1988. Our main results focus on an optimal control problem associated with this model, where the short frequency, understood as a bounded time-dependent function, is treated as the control variable. The cost function accounts for a person's pain at a given final time and incorporates additional criteria. We demonstrate the uniqueness of the optimal control and establish the bang-bang nature of the control. In a previous section, we extend the mathematical analysis of the model developed by Britton and Skevington by presenting a series of mathematical inequalities. These inequalities strengthen the model's alignment with the principal requirements for reproducing the core structure of Pain Theory.
\end{abstract}


\section{Introduction and some preliminaries}
\label{sec:Intro}
The purpose of this paper is to analyze an optimal control problem associated with the so-called "Pain Gate Theory" proposed in 1965 by R. Melzack and P. Wall (\cite{MW1}) and later confirmed and extended in their 1982 revision (see \cite{MW2}). A mathematical model related to this theory was proposed by N.F. Britton and S.M. Skevington \cite{BS} (see also the review\cite{Na}). This theory forms the basis of the so-called "neurostimulators," devices implanted in patients with chronic pain to alleviate their discomfort. Our mathematical study consists of two parts. In the first part, we demonstrate that, under appropriate structural assumptions, the corresponding coupled system has a unique solution. Additionally, we analyze various qualitative properties of the solutions, focusing particularly on their continuous and monotonic dependence on the data. In the second part, we examine the optimal control problem associated with the equation below for the potentials $\bV(t) = \big (\rV_{\rE}(t), \rV_{\rI}(t), \rV_{\rT}(t)\big )^{\tt t}$, considering the short frequency $x_{s}$ as the control variable. The system under consideration can be expressed as follows:
\begin{equation}
\left \{
\begin{array}{l}
\derp{\bV} (t)=\bff \big (\bV(t),x_{s}(t)\big ),\quad 0\le t,\\ [.175cm]
\bV(0)=\bV_{0},
\end{array}
\right .
\label{eq:CDS}
\end{equation}
where $ \bV_{0}=(\rV_{\rE 0},\rV_{\rI 0},\rV_{\rT 0})^{\tt t}$ and
\begin{equation}
\bff \big (\bV,x_{s}\big )=
\left (
\begin{array}{l}
\big (-\rV_\rE+\rV_{\rE 0} + g_{s\rE}(x_s)\big )\tau_\rE ^{-1} \\ [.1cm]
\big (-\rV_\rI+\rV_{\rI 0}+\rG_{\rT}(\rV_{\rT}+\rC_{\rI}\big )\tau_\rI^{-1} \\ [.1cm]
\big (-\rV_\rT+\rV_{\rT 0}-\rG_{\rI}(\rV_{\rI})+\rG_{\rE}\big (\rV_{\rE}\big ) + g_{s\rT}(x_s) + g_{l\rT}(x_l)\tau_\rT^{-1}
\end{array}
\right )
\label{eq:filed1}
\end{equation}
(see more details in \eqref{eq:Gfunctions} and \eqref{eq:constantsterms} below).
The formulation as an optimal control problem arises from minimizing a functional in which the pain at the final time is assumed to be an increasing function of the $\rV_{\rT}(t_{f})$ potential (for a given final time $t_{f}>0)$)
\begin{equation}
\rW(\bV_{0})\doteq \min_{x_{s}\in \cU}\J_{\bV_{0}} (x_{s}),\quad \J_{\bV_{0}} (x_{s})\doteq \int^{t_{f}}_{0}\left (
\dfrac{1}{2}|\rV_{\rT}(\sigma)|^{2}- \rQ g_{s\rT}\big (x_{s}(\sigma)\right )d\sigma+\rS\big (\rV_{\rT}(t_{f})\big).
\label{eq:minimization}
\end{equation}
Here, the set of controls is given by $\cU=\rL^{\infty}\big (\R_{+}:[\underline{x_{s}},\overline{x_{s}}]\big )$, $\rQ$ is a given positive constant, and $\rS$ is a $\cC^{1}$ function satisfying {\em $0<\rS_{-}'\le \rS'\big (\rV_{\rT}\big )\le \rS'_{+}$}. From a practical point of view, and as is common in Control Theory, the first term of the functional $\J_{\rV}$ represents the "cost of the state" (which, indirectly, implies that we aim to minimize the electrical energy). The second term represents a quantitative measure of pain at a prescribed (but arbitrary) fixed final time.

This functional is not convex with respect to the control due to the nonlinearity of the state equations; thus, the uniqueness of the minimum is not guaranteed by the abstract theory. Nevertheless, we will prove that the Pontryagin Principle applies. Consequently, we will demonstrate the existence of a "bang-bang" type control and explicitly characterize the unique control that has a single point of discontinuity.

The organization of this paper is as follows: Section 2 is devoted to a description of the modeling, leading to the coupled system under consideration. The well-posedness of such a system and, especially, the monotonicity properties of the components of the solution with respect to the data are presented in Section 3, where we provide several quantitative estimates that give a sharper formulation of some of the results in \cite{BS}. We also present a series of figures showing the simulation of some special cases using {\tt MATLAB}. Finally, Section 4, containing the main results of this paper, is devoted to the study of the optimal control problem.

We mention that our mathematical approach to this model contains several improvements over the important paper by \cite{BS}, and, for instance, includes a series of figures resulting from some {\tt MATLAB} approximations. Moreover, the consideration of the associated control problem is entirely new compared to the previous literature on the model. Clearly, the consideration of the control study of the model leads to a very useful  methodology in this process since it can be adapted to each patient in a more singularized way (physical characteristics, pain sensitivity, etc.).
\section{On the mathematical model}
\label{sec:model}
The key ideas for the construction of a mathematical model related to pain were presented in the seminal paper by N.F. Britton and S.M. Skevington in 1988 \cite{BS}. As stated in that work (where some earlier references are cited), any theory of pain must be able to account for the following observations:
(i) Increased stimulation of the small nerve fibers in the skin usually increases pain.
(ii) Increased stimulation of the large nerve fibers may increase pain temporarily, but in the long term, it may relieve it.
(iii) Pain relief can be achieved through electrical stimulation of the grey matter in the midbrain.
(iv) Injuries that would normally cause severe pain sometimes cause little or no pain, or the onset of pain is delayed.
(v) In some cases, the anticipation of pain is enough to raise anxiety levels, thereby intensifying the perception of pain.

\begin{figure}
	\centering
	\includegraphics[width=12cm]{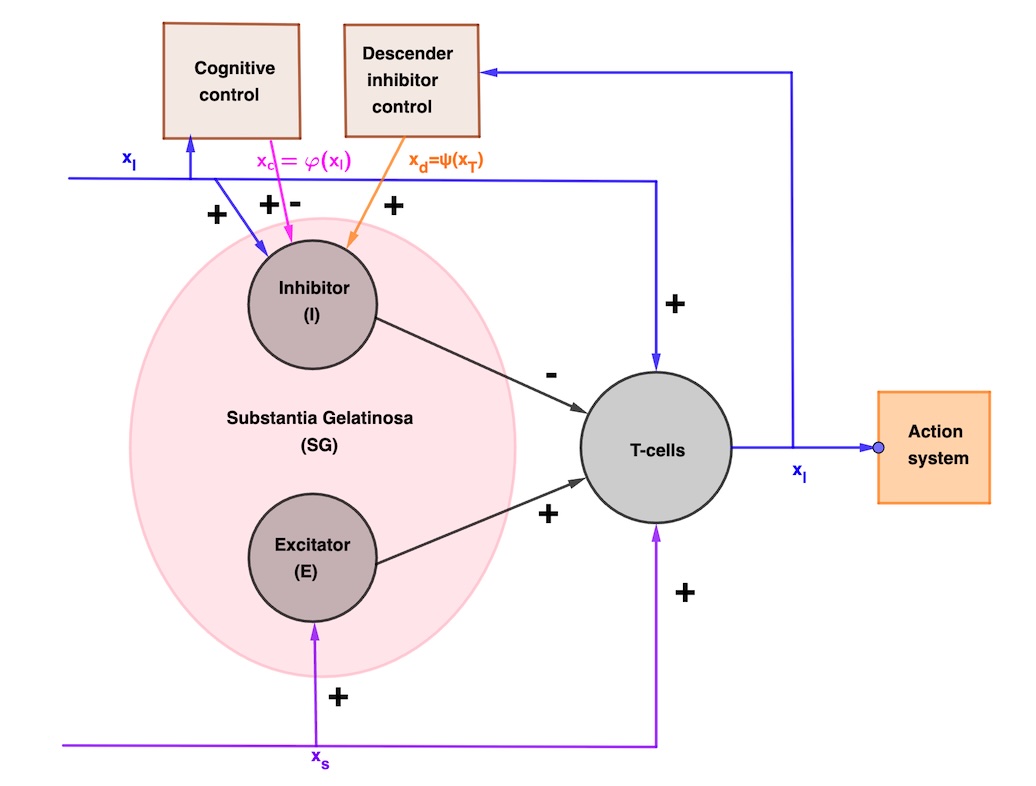}\\ [-.5cm]
	\caption{\small Basic  design of the pain gate}
	\label{fig:dolor2}
\end{figure}

The mathematical model used in this paper is built on the above ideas and can be sketched as indicated in Figure \ref{fig:dolor2} (adaptation from \cite{MW2} and \cite{BS}):  firing frequency $x$ in the pathways due to a firing frequency $x_s$, in the small fibers and  $x_l$ in the large fibers of the considered area of skin.  We will assume that the frequency of the outputs 
from the  cognitive control and the descending inhibitory control are strictly increasing functions of the inputs, i.e.,
\begin{equation} 
\overbrace{x_d = \varphi(x_\rT),}^{\hbox{\tiny  {descending inhibitory control}}} \qquad  \overbrace{x_c = \psi(x_l)}^{\hbox{\tiny  cognitive control}},   
\label{eq:frecuencies} 
\end{equation}
with $\varphi$, $\psi$  strictly increasing such that $\varphi(0) = 0$, $\psi(0) = 0$.
\par
Following the ideas of Wilson \& Cowan \cite{WC}, we will assume that each T-cell is stimulated by one large and one small afferent nerve fiber from the skin, as well as one inhibitory and one excitatory SG cell (substantia gelatinosa). We will use the potentials of the T-cell, $\rV_\rT$, and those of the inhibitory and excitatory SG cells, $\rV_\rI$ and $\rV_\rE$, respectively. The frequencies $x_\rI$ and $x_\rT$ at which these cells fire are functions of the slow potentials. 
\begin{equation} 
x_\rT = f_\rT(\rV_\rT), \qquad x_\rI = f_\rI(\rV_\rI), \qquad x_\rE = f_\rE(\rV_\rE).
\label{eq:frecuenciaspotenciales3}
\end{equation}
We assume that these functions vanish for suitable values of $\rV$ smaller than a certain threshold $\rV_{u}$ and are strictly increasing for values larger than the threshold.
$$
f_{k}(\rV)=-\rL\big (\rV-\rV_{uk}\big )_{+},\quad \rV\in\R \quad \big (r_{+}\doteq \max\{r,0\}\big ),
$$
(see Remark \ref{obser:ejemplo} below).
\par
The potentials $\rV$ are assumed to depend on the frequencies of impulses arriving at their dendrites from various sources. We will also assume that the properties of the dendrites are constant over time.
According to An der Heiden \cite{An}, the effect of an input frequency $x_j$ to an excitatory or inhibitory synapse of a cell with potential $\rV_k$ will be to raise it by
\begin{equation} 
\Phi_{jk} = \alpha_{jk} \int^t_{-\infty}{h_{jk}(t-\tau) g[x_j(\tau)] d\tau},
\label{eq:Phi}
\end{equation}
with $\alpha_{jk} = 1$ for an excitatory and $-1$ for an inhibitory synapse. Here, $h_{jk}$ is a positive monotone decreasing function, and $g_{jk}$ is a bounded strictly monotone increasing function such that $g_{jk}(0)=~0$.
A simple choice for $h_{jk}$ corresponds to when it is associated with an RC-network composed of resistors and capacitors.
\begin{equation} 
\rC_{k}\derp{h_{jk}}(t)+\dfrac{h_{jk}(t)}{\rR_{k}}=0\quad \Rightarrow\quad 	h_{jk}(t) =h_{jk}(0)  \text{ exp} \left( -\frac{t}{\rR_{k}\rC_{k}} \right).
\label{eq:h}
\end{equation}
In that case, the total input effect on the cell $k$ is of the form
\begin{equation}
\rV_k = \overbrace{\rV_{k0}}^{\hbox{\tiny initial potential}} + \overbrace{\sum_j \Phi_{jk}}^{\hbox{\tiny synapse effect}},
\label{eq:potenciales}
\end{equation}
where we have iterated over the inputs $j$ arriving at or departing from the cell $k$.

From	\eqref{eq:Phi} and \eqref{eq:h}, differentiating \eqref{eq:potenciales} we get
\begin{equation}
\derp{\rV}_k =\bigg [\sum_j \alpha_{jk} h_{jk}(0)g_{jk}(x_j)-\dfrac{1}{\rR_{k}\rC_{k}}\overbrace{\sum_j \Phi_{jk}}^{\rV_{k}-\rV_{k0}}\bigg ].
\label{eq:variaciongenericapotenciales}
\end{equation}
\par
\vspace*{-.2cm}
Taking $h_{jk}(0) = \dfrac{1}{\tau_{k}}$, with $\tau_{k} = \rR_{k}\rC_{k}$ (a membrane constant), the dynamics of the potentials are given by
\begin{equation}
\left \{
\begin{array}{ll}
\tau_\rE \derp{\rV}_\rE =&\hspace*{-.2cm} -(\rV_\rE-\rV_{\rE 0}) + g_{s\rE}(x_s), 
\\ [.05cm]
\tau_\rI \derp{\rV}_\rI =&\hspace*{-.2cm} -(\rV_\rI-\rV_{\rI 0}) +g_{l\rI}(x_l)  + g_{d\rI}(x_d)+ \alpha_{c\rI}g_{c\rI}(x_c) ,
\\ [.05cm]
\tau_\rT \derp{\rV}_\rT =&\hspace*{-.2cm} -(\rV_\rT-\rV_{\rT 0}) + g_{s\rT}(x_s) + g_{l\rT}(x_l) + g_{\rE\rT}(x_\rE) - g_{\rI\rT}(x_\rI), 
\end{array}
\right .
\label{eq:systemprevious} 
\end{equation}
(see Figure \ref{fig:dolor2}).
Here, $\alpha_{c\rI} \in [-1,1]$ is a cognitive control, which is positive for an excitatory input, negative for an inhibitory input, and zero for no input from the cognitive control. For a deeper understanding of the coefficients' meaning (and, in particular, why alpha is associated with the intensity of cognitive control), we send the reader to  \cite{BCS} and its references. 
\par
\vspace*{.1cm}
According to \eqref{eq:frecuencies} and \eqref{eq:frecuenciaspotenciales3}, we obtain the nonlinear system modeling the dynamics of the potentials $\bV = (\rV_{\rE}, \rV_{\rI}, \rV_{\rT})^{\tt t}$ in terms of the known inputs $x_s$ and $x_l$
\begin{equation}
\left \{
\begin{array}{ll}
	\tau_\rE \derp{\rV}_{\rE} =&\hspace*{-.2cm} -(\rV_{\rE}-\rV_{\rE 0}) + g_{s{\rE}}(x_s), \\ 
\tau_\rI \derp{\rV}_{\rI} =&\hspace*{-.2cm} -(\rV_{\rI}-\rV_{\rI 0}) + g_{l\rI}(x_l) + g_{d\rI}\big (\varphi\big [f_{\rT}(\rV_{\rT})\big ]\big )+\alpha_{c\rI}g_{c{\rI}}\big [\psi(x_l)\big ] ,\\ 
\tau_\rT \derp{\rV}_{\rT} =& \hspace*{-.2cm}-(\rV_{\rT}-\rV_{\rT 0}) + g_{s\rT}(x_s) + g_{l\rT}(x_l) + g_{\rE\rT}\big [f_\rE(\rV_\rE)\big ] - g_{\rI\rT}[f_\rI(\rV_\rI)],
\end{array}
\right.
\label{eq:system3} 
\end{equation}
assuming that we prescribe an initial behaviour $\bV(0)=\big (\rV_{\rE}^{0},\rV_{\rI}^{0},\rV_{\rT}^{0}\big )^{\tt t}\in\R^{3}$. In particular, the equation
\begin{equation}
\tau_\rT \derp{\rV}_{\rT} = -(\rV_{\rT}-\rV_{\rT 0}) + g_{s\rT}(x_s) + g_{l\rT}(x_l) + g_{\rE\rT}\big [f_\rE(\rV_\rE)\big ] - g_{\rI\rT}[f_\rI(\rV_\rI)],
\label{eq:ecuacionVt3} 
\end{equation}
collects the action of the frequencies $x_{s}$ and $x_{l}$, as well as the potentials $\rV_{\rI}$ and $\rV_{\rE}$, of the inhibitory and excitatory SG cells on the variation of the potential $\rV_{\rT}$ of the T cell (see the illustrations in Figures \ref{fig:VariacionVtfrecuenciacorta} and \ref{fig:VariacionVtmodo} below).

\begin{rem}\rm In 1996, N.F. Britton, M.A.J. Chaplain, and S.M. Skevington \cite{BCS} added to the "Gate Control Theory of Pain" a variant that takes into account the so-called wind-up mechanism, in which certain receptors (the N-methyl-D-aspartate, NMDA), which are highly relevant in pain sensitivity, reactivate the process (wind-up). Along with \eqref{eq:system3}, they introduced a new unknown, $\rV_{m}$, modeling the mid-brain input potential, and a new frequency, $x_{m} = f_{m}(\rV_{m})$. The new system becomes 
\begin{equation}
\left \{
\begin{array}{ll}
\tau_\rE \derp{\rV}_{\rE} =&\hspace*{-.2cm} -(\rV_{\rE}-\rV_{\rE 0}) + g_{s{\rE}}(x_s,\rV_{\rE}), \\ 
\tau_m \derp{\rV}_{m} =& \hspace*{-.2cm}-(\rV_{m}-\rV_{m 0}) + g_{m\rT}(x_\rT), \\
\tau_\rI \derp{\rV}_{\rI} =&\hspace*{-.2cm} -(\rV_{\rI}-\rV_{\rI 0}) + g_{l\rI}(x_l) + g_{m\rI}(x_{m}),\\ 
\tau_\rT \derp{\rV}_{\rT} =& \hspace*{-.2cm}-(\rV_{\rT}-\rV_{\rT 0}) +  g_{s\rT}(x_s)+g_{l\rT}(x_l) + g_{\rE\rT}(x_\rE) - g_{l\rT}(x_{\rI})-g_{m\rT}(x_{m}).\\
\end{array}
\right.
\label{eq:system4} 
\end{equation}
By using the frequency potential's dependence 
\begin{equation} 
x_\rT = f_\rT(\rV_\rT), \qquad x_\rI = f_\rI(\rV_\rI), \qquad x_\rE = f_\rE(\rV_\rE), \qquad x_m = f_m(\rV_m),
\label{eq:frecuenciaspotenciales4}
\end{equation}
we get the nonlinear system modeling the dynamics of the potentials $\bV= \big(\rV_{\rE}, \rV_{m}, \rV_{\rI}, \rV_{\rT}\big)^{\tt t}$ in terms of the known frequencies $x_s$ and $x_l$, assuming the initial behaviour $\bV(0) = \big(\rV_{\rE}^{0}, \rV_{m}^{0}, \rV_{\rI}^{0}, \rV_{\rT}^{0}\big)^{\tt t}$, similarly to what was done for \eqref{eq:system3}. Notice that now the equation
\begin{equation}
\tau_\rT \derp{\rV}_{\rT} = -(\rV_{\rT}-\rV_{\rT 0}) + g_{s\rT}(x_s) + g_{l\rT}(x_l) + g_{\rE\rT}\big [f_\rE(\rV_\rE)\big ] - g_{\rI\rT}[f_\rI(\rV_\rI)]- g_{m\rT}[f_m(\rV_m)],
\label{eq:ecuacionVt4} 
\end{equation}
collects the action of $x_{s}$ and $x_{l}$, as well as $\rV_{\rI}$, $\rV_{\rE}$, and $\rV_{m}$, on the variation of the potential $\rV_{\rT}$ of the T cell.
The qualitative treatment of both systems is quite similar. In our presentation, we have chosen system \eqref{eq:system3} since it is easier to understand, but very slight changes lead to similar results for the other system. $\fin$

\label{obser:cuatropotenciales}
\end{rem}

\section{Dynamic of the states of the control problem}
\label{sec:dynamics}
In this section, we will study the dynamics of the states $\bV = \big (\rV_{\rE}, \rV_{\rI}, \rV_{\rT}\big )^{\tt t}$ of the system \eqref{eq:system3}, assuming the controls of the problem are known. We can reformulate the system in vector form as
\begin{equation}
\derp{\bV} (t) =\bff (\bV(t))
\label{eq:sistemadinamico}
\end{equation}
with
\begin{equation}
\bff (\bV)=
\left (
\begin{array}{l}
\big (-(\rV_\rE-\rV_{\rE 0}) + g_{s\rE}(x_s)\big )\tau_\rE ^{-1} \\ [.1cm]
\big (-(\rV_\rI-\rV_{\rI0}) + g_{l\rI}(x_l) + \rG_{\rT}(\rV_\rT\big ) + \alpha_{c\rI}g_{c\rI}\big [\psi(x_l)\big ]\big )\tau_\rI^{-1} \\ [.1cm]
\big (-(\rV_{\rT}-\rV_{\rT 0}) + g_{s\rT}(x_s) + g_{l\rT}(x_l) + \rG_{\rE\rT}(\rV_\rE) - \rG_{\rI}(\rV_\rI)\big )\tau_\rT^{-1}
\end{array}
\right )
\label{eq:field} 
\end{equation}
by assuming the given strictly increasing and bounded functions
\begin{equation}
\rG_{\rT}(\rV_{\rT})\doteq g_{d\rI}\big (\varphi\big [f_\rT(\rV_{\rT})\big ]\big ),~\rG_{\rE}(\rV_{\rE})\doteq g_{\rE\rT}\big [f_\rE(\rV_{\rE})\big ] \quad \hbox{and}\quad \rG_{\rI}(\rV_{\rI})\doteq g_{\rI\rT}[f_\rI(\rV_{\rI})].
\label{eq:Gfunctions}
\end{equation}
Such functions depend on the choice of $\alpha_{c\rI} \in [-1,1]$ and the frequencies $x_s$ and $x_l$. We will make the following assumptions.
\begin{itemize}
\item[(F)]  $f_k\in \cC(\R,\R_+)\cap \cC^{1}((\rV_{uk},\infty):\R_+)$, with 
$f_k(s) = 0$, for $s < \rV_{uk}$, for some threshold $\rV_{uk}$, and $f_k(s)$ is strictly increasing for $s > \rV_{uk}$. \item[(G)] For any $j, k$, the functions $g_{jk} \in \cC^1(\R_{+}, \R_{+})$ are bounded and strictly increasing such that $g_{jk}(0) = 0$. \item[(H)] The functions $\varphi, \psi \in \cC^1(\R, \R_+)$ are strictly increasing and such that $\varphi(0) = \psi(0) = 0$. 
\item[(I)] $x_{s} \in \rL^{1}(0, t_{f}; \R)$. 
\end{itemize}

The following result extends the similar version obtained in \cite{BS} to the case of time-dependent small frequencies.
\begin{theo}[Existence and uniqueness of solutions] Under conditions {\rm (F)}, {\rm (G)}, {\rm (H)} and {\rm (I)} there exists a unique absolutely continuous solution on the interval $\left[ 0,t_{f}\right]$, for any arbitrarily fixed $t_{f} > 0$, to the initial value problem. Moreover, this solution is bounded for any $t \in \left[ 0,t_{f}\right]$.
\label{teo:existenciaunicidad}
\end{theo}
\proof Function  $\bff(\bV)$ is continuous and all the terms of the gradient matrix 
\begin{equation}
\rD_{\bV}\bff (\bV)=
\left (
\begin{array}{c|cc}
-\tau_\rE ^{-1} & & \\  [.1cm] \hline 
& -\tau_{\rI}^{-1} & g_{d\rI}'\big (\varphi'\big [f_t'(\rV_\rT)\big ]\big )\tau_{\rT}^{-1}  \\ [.1cm]
g_{\rE\rT}'[f_\rT '(\rV_\rE)]\big ) & - g_{\rI\rT}'[f_\rI '(\rV_\rI)]\big ) &-\tau_\rT^{-1}
\end{array}
\right )
\label{eq:gradientecampos} 
\end{equation}
are bounded, so the function $\bff(\bV)$ is, in fact, globally Lipschitz continuous. Thus, the existence and uniqueness of solutions follow from the Picard-Lindelöf Theorem, even if the frequency data are discontinuous but belong to some suitable space, such as $x_{s} \in \rL^{1}(0,t_{f}; \R)$ (see, e.g., \cite[Theorem 3.2.1]{BP}). On the other hand, it is straightforward to obtain the following inequality estimates 
$$
\big |\tau_{k}\derp{\rV}_{k}(t)+\rV_{k}(t)\big |\le \rM_{k}\quad \Rightarrow \quad |\rV_{k}(t)|\le \big |\rV_{k}(0)-\rM_{k}\big |e^{-\frac{t}{\tau_{k}}}+\rM_{k}\le \big |\rV_{k}(0)-\rM_{k}\big |+\rM_{k}
$$ 
for some positive constant $\rM_{k}$, which depends only on the data. $\fin$
\par
\begin{rem}\rm \label{Rm:Special case} Following the reference \cite{BCS}, it is relevant to illustrate the results for the choice of
$$
f_{k}(x)=-\dfrac{\rK}{\rV_{k0}}\big (x-\rV_{uk}\big )_{+},\quad x\in\R,
$$
where $\rK$ is a positive constant and $\rV_{uk}$ is the pain threshold associated with the potential $\rV_{k}$. For the sake of presentation, in the next examples, we will always take $\rV_{u} = -55$. We will also take $g_{jk}(x) = c_{jk} \tanh(x)$, with $c_{jk} \in \R_{+}$ (see \cite{BCS}). In addition, in all the following examples, we replace assumption {\rm (I)} with the simpler case
\begin{equation}
x_{s}\text{ is constant in time.}  
\label{eq:constant in time}
\end{equation}
Finally, we will consider the simple cases $\varphi(x)=x^{n}$ and $ \psi(x)=x^{m}$ for some odd exponents $n$ and~$m.\fin$ 
\label{obser:ejemplo}
\end{rem}
Notice that from  \eqref{eq:field}  we get that the excitatory potential is uncoupled:
$$
\tau_\rE \dot{\rV}_\rE =-(\rV_\rE-\rV_{\rE0}) + g_{s\rE}(x_s). 
$$
More precisely, assuming \eqref{eq:constant in time} we obtain
\begin{equation}
\rV_{\rE}(t)=\big (\rV_{\rE}(0)-\rV_{\rE}(\infty)\big )e^{-\frac{t}{\tau_{\rE}}}+\rV_{\rE}(\infty),\quad t\ge 0, 
\label{eq:Ve}
\end{equation}
with 
\begin{equation}
\rV_{\rE}(\infty)\doteq  \rV_{\rE 0} + g_{s\rE}(x_s), 
\label{eq:Veinfinito}
\end{equation}
In other words, the stationary excitatory potential depends solely on the corresponding initial potential and the short frequency. 
\par
\medskip
An example (following the special case indicated in Remark \ref{Rm:Special case}) of the qualitative behavior of solutions derived from the above-mentioned a priori estimates is shown in Figure \ref{fig:graficaspotenciales} below.
\begin{figure}[ht]
\centering
\includegraphics[width=11cm]{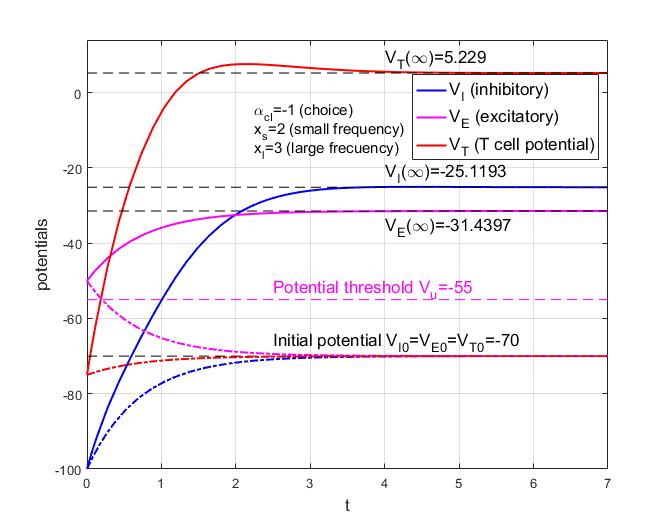} \\ [-.35cm]
\caption{Qualitative representation of the potentials $\rV_{\rI},\rV_{\rE}$ y $\rV_{\rT}$}
\label{fig:graficaspotenciales}
\end{figure}
The following series of technical results are generalizations of similar statements made in \cite{BS}, but expressed only in qualitative terms:
\begin{theo}
Assume \eqref{eq:constant in time}. Then the dynamical system associated with \eqref{eq:field} has a unique stationary state $\bV(\infty) = \big(\rV_{\rI}(\infty), \rV_{\rE}(\infty), \rV_{\rT}(\infty)\big)$. Moreover, the property $\rV_{uk} \le \rV_{k}(\infty)$ holds for all the potential components.
\label{theo:phaseplain}
\end{theo}
\proof Once we obtain the result for $\rV_{\rE}(\infty)$ (see \eqref{eq:Veinfinito}), the rest of the components of the stationary potential $\big(\rV_{\rI}(\infty), \rV_{\rE}(\infty)\big)$, with $\rV_{\rI}(\infty) > \rV_{u\rI}$ and $\rV_{\rT}(\infty) > \rV_{u\rT}$, are solutions of the system
$$
\left \{
\begin{array}{l}
-(\rV_\rI-\rV_{\rI 0}) + g_{l\rI}(x_l) +\rG_{\rT}\big (\rV_\rT\big ) + \alpha_{c\rI}g_{c\rI}\big [\psi(x_l)\big ]=0,\\ [.15cm]
-(\rV_\rT-\rV_{\rT 0}) + g_{s\rT}(x_s) + g_{l\rT}(x_l) + \rG_{\rE}\big (\rV_\rE(\infty)\big ) - \rG_{\rI}\big (\rV_\rI\big )=0.
\end{array}
\right.
$$
We can rewrite such a system as
\begin{equation}
\left \{
\begin{array}{l}
\rV_{\rI}-\rG_{\rT}(\rV_{\rT})=\rC_{\rI},\\ [.15cm]
\rG_{\rI}(\rV_{\rI})+\rV_{\rT}=\rC_{\rT}+\rG_{\rE}\big (\rV_{\rE}(\infty)\big ),
\end{array}
\right.
\label{eq:sistemaestacionarionotacional}
\end{equation}
where  the constant terms are given by
\begin{equation}
\left \{
\begin{array}{ll}
\rC_{\rI}\doteq \rV_{\rI 0} + g_{l\rI}(x_l) + \alpha_{c\rI}g_{ci}\big [\psi(x_l)\big ]&\quad \hbox{( independent of the short frequency $x_{s}$)}, \\ [.15cm]
\rC_{\rT}\doteq \rV_{\rT 0} + g_{s\rT}(x_s) + g_{l\rT}(x_l)&\quad \hbox{(independent of the choice of $\alpha_{c\rI}$)},
\end{array}
\right.
\label{eq:constantsterms}
\end{equation}
(see also \eqref{eq:Gfunctions}). It is easy to see that from the system \eqref{eq:sistemaestacionarionotacional} we deduce that the following scalar equation on $\rV_{\rT}$
\begin{equation}
\rV_{\rT}+\rG_{\rI}\big (\rG_{\rT}(\rV_{\rT})+\rC_{\rI}\big )=\rC_{\rT}+\rG_{\rE}\big (\rV_{\rE}(\infty)\big )
\label{eq:equlibriumexistence}
\end{equation}
has a unique solution $\rV_{\rT}(\infty)$ since the mapping $\rV_{\rT} \mapsto \rV_{\rT} + \rG_{\rI} \big( \rG_{\rT}(\rV_{\rT}) + \rC_{\rI} \big)$ is a strictly increasing function, and its range is the entire space. Then we get $\rV_{\rI}(\infty)=\rC_{\rI}+\rG_{\rT}(\rV_{\rT}(\infty)).\fin$ 
\par
\medskip
Once again, we can illustrate the above result (following the special case indicated in Remark~\ref{Rm:Special case}) as shown in Figure \ref{fig:graficaspotenciales}.

\begin{figure}[ht]
\centering
\includegraphics[width=10.5cm]{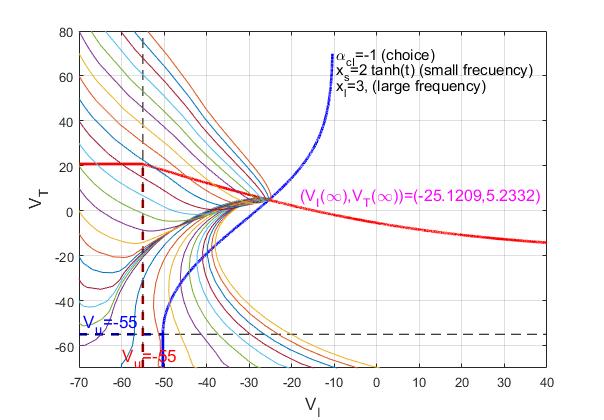} \\ [-.35cm]
\caption{Phase plain for the potentials $\rV_{\rI}$ and $\rV_{\rT}$}
\label{fig:planofases}
\end{figure}

The following result was obtained in \cite{BS}. 
\begin{theo}[\cite{BS}] Assume \eqref{eq:constant in time}. Then, the stationary state $\bV(\infty)$ is asymptotically stable. 
\end{theo}
\proof From \eqref{eq:field} we consider the linearized system
$$
\derp{\bV}(t)=\rD_{\bV}\bff \big (\bV(\infty)\big )\bV(t)
$$
(see \eqref{eq:gradientecampos}). The corresponding characteristic polynomial is 
$$
\begin{array}{ll}
\cP(\lambda)&\hspace*{-.2cm}=\big |\rD_{\bV}\bff \big (\bV(\infty)\big )-\lambda \rI\big |\\  [.2cm]
&\hspace*{-.2cm}=
\left |
\begin{array}{c|cc}
-\tau_\rE ^{-1}-\lambda  &  & \\  [.1cm] \hline
& -\tau_{\rI}^{-1}-\lambda  & \rG'_{\rI}\big (\rV_{\rI}(\infty)\big )\tau_{\rI}^{-1}  \\ [.1cm] 
- \rG_{\rE}'\big (\rV_{\rE}(\infty)\big )\tau_{\rT}^{-1} &
- \rG_{\rT}'\big (\rV_{\rT}(\infty)\big )\tau_{\rT}^{-1} & -\tau_\rT^{-1}-\lambda 
\end{array}
\right |\\ [.9cm]
&
=-\big (\lambda +\tau_{\rE}^{-1}\big )\left [\big (\lambda +\tau _{\rI}^{-1}\big )\big (\lambda +\tau_{\rT}^{-1}\big )+  \rG'_{\rI}\big (\rV_{\rI}(\infty)\big )\rG_{\rT}'\big (\rV_{\rT}(\infty)\big )\tau_{\rT}^{-1}\tau_{\rI}^{-1}\right ] \\ [.2cm]
&
=-\big (\lambda +\tau_{\rE}^{-1}\big )\left [\lambda ^{2}+\big (\tau_{\rI}^{-1}+\tau_{\rT}^{-1}\big )\lambda +\tau_{\rI}^{-1}\tau_{\rT}^{-1}+
\rG'_{\rI}\big (\rV_{\rI}(\infty)\big )\rG_{\rT}'\big (\rV_{\rT}(\infty)\big )\tau_{\rT}^{-1}\tau_{\rI}^{-1}\right ].
\end{array}
$$
This polynomial has a real root $\lambda_{1} = -\tau_{\rE}^{-1}$ and two other roots, $\lambda_{2}$ and $\lambda_{3}$, with negative real parts. Then, by the Grobman-Hartman Theorem (see e.g. \cite{Am}), we obtain the result. $\fin$
\begin{rem}\rm Assuming \eqref{eq:constant in time}, by applying the arguments from \cite{CD}, it is possible to prove the following asymptotic estimate:
$$
\limsup_{t\rightarrow \infty} e^{\min\{{\rm Re}(\lambda_{1}),{\rm Re}(\lambda_{2}),{\rm Re}(\lambda_{2})\} t}\n{\bV(t)-\bV(\infty)}\in \R_{+}.
$$	
$\fin$
\end{rem}

\par
In the second part of this section, we will study the dependence of the potentials on the parameters. The following technical result is an easy consequence of the argument used in the proof of Theorem \ref{theo:phaseplain}.

\begin{lemma} Assume \eqref{eq:constant in time} and consider the identities
$$
\left \{
\begin{array}{l}
\rV_{\rI}-\rG_{\rT}(\rV_{\rT})=\rC_{\rI},\\ [.15cm]
\rG_{\rI}(\rV_{\rI})+\rV_{\rT}=\rC_{\rT}+\rG_{\rE}(\rV_{\rE}),
\end{array}
\right .
\quad \hbox{and} \quad 
\left \{
\begin{array}{l}
\rVh_{\rI}-\rG_{\rT}(\rVh_{\rT})=\rCh_{\rI},\\ [.15cm]
\rG_{\rI}(\rVh_{\rI})+\rVh_{\rT}=\rCh_{\rT}+\rG_{\rE}(\rVh_{\rE}),
\end{array}
\right .
$$
with functions $\rG_{\rI}$ and $\rG_{\rT}$ defined in \eqref{eq:Gfunctions}. Then, there exist intermediate values $\theta_{\rT}$, between $\rV_{\rT}$ and $\rVh_{\rT}$, and $\theta_{\rI}$, between $\rV_{\rI}$ and $\rVh_{\rI}$, such that the following expressions hold:
\begin{equation}
\rV_{\rI}-\rVh_{\rI}=
\dfrac{\big (\rC_{\rI}-\rCh_{\rI}\big )-\rG_{\rT}'(\theta_{\rT})\big (\rC_{\rT}+\rG_{\rE}(\rV_{\rE})-\rCh_{\rT}-\rG_{\rE}(\rVh_{\rE})\big )}{1+\rG_{\rI}'(\theta_{\rI})\rG_{\rT}'(\theta_{\rT})}
\label{eq:variacionVi}
\end{equation}
and 
\begin{equation}
\rV_{\rT}-\rVh_{\rT}=
\dfrac{\big (\rC_{\rT}+\rG_{\rE}(\rV_{\rE})-\rCh_{\rT}-\rG_{\rE}(\rVh_{\rE})\big )-\rG_{\rI}'(\theta_{\rI})\big (\rC_{\rI}-\rCh_{\rI}\big )}{1+
	\rG_{\rI}'(\theta_{\rI})\rG_{\rT}'(\theta_{\rT})}.
\label{eq:variacionVt}
\end{equation}
$\fin$
\label{lema:variacionVtVi}
\end{lemma}
\begin{rem}\rm
We point out that the expressions $\rG_{\rI}'(\theta_{\rI})$ and $\rG_{\rT}'(\theta_{\rT})$ are well-defined since the stationary potentials are beyond the corresponding threshold value. $\fin$
\end{rem}
The following result allows us to obtain a quantitative indication of the monotone dependence of~$\rV_{\rT}(\infty)$.
\begin{prop}[Monotone dependence for $\rV_{\rT}(\infty)$] \
\par
\noindent
Assume \eqref{eq:constant in time}. Then:
\par
\noindent
a) {\sc Dependence with respect to the short frequency $x_{s}$}. Given $\rV_{\rE 0}$, $x_{l}$, and $\alpha_{c\rI}$, the stationary potential $\rV_{\rT}(\infty;x_{s})$ is a monotone function of $x_{s}$ with the same monotonicity as $g_{s\rT}(x_{s})$. Moreover, we have the inequality 
\begin{equation}
\big [\rV_{\rT}\big (\infty;x_{s}\big )-\rV_{\rT}\big (\infty;\xh_{s}\big )\big ]_{+}\le 
\big [g_{s\rT}(x_{s})+\rG_{\rE}\big (\rV_{\rE 0}+g_{s\rE}(x_{s})\big )-g_{s\rT}(\xh_{s})-\rG_{\rE}\big (\rV_{\rE 0}+g_{s\rE}(\xh_{s})\big )\big ]_{+}.
\label{eq:Vtxs}
\end{equation}
b) {\sc Dependence with respect to $\alpha_{c\rI}$}. Given $\rV_{\rE 0}$, $x_{s}$ and $x_{l}$,  the stationary potential $\rV_{\rT}(\infty;\alpha_{c\rI})$ is a decreasing function on $\alpha_{c\rI}$. Moreover, the following inequality holds
\begin{equation}
\big [\rV_{\rT}\big (\infty;\alpha_{1}\big )-\rV_{\rT}\big (\infty;\alpha_{2}\big )\big ]_{-}\le\left ( \sup_{\theta}\dfrac{\rG_{\rI}'(\theta)}
{1+ \rG_{\rI}'(\theta)\rG_{\rT}'(\theta)}\right )g_{c\rI}\big [\psi (x_{l})\big ]\big [\alpha_{1}-\alpha_{2}\big ]_{+}).
\label{eq:Vtaci}
\end{equation}
Here we used the notation  $r_{+}=\max\{r,0\}$ and $r_{-}=\max\{-r,0\}.$
\label{prop:variationVtparameters}
\end{prop} 
\proof
\par
\noindent
a)  Given two values $x_{s}$ and $\xh_{s}$, we consider the associated potentials $\rV_{\rI}\big (\infty;x_{s}\big )$, $
\rV_{\rI}\big (\infty;\xh_{s}\big )$ y
$\rV_{\rT}\big (\infty;x_{s}\big )$, $\rV_{\rT}\big (\infty;\xh_{s}\big )$, as well as  the constant terms $\rC_{\rI}(x_{s}), \rC_{\rI}(\xh_{s})$ and 
$\rC_{\rT}(x_{s}), \rC_{\rT}(\xh_{s})$, given in the corresponding system 
\eqref{eq:sistemaestacionarionotacional}. Then \eqref{eq:variacionVt} leads to
$$
\begin{array}{ll}
	\rV_{\rT}\big (\infty;x_{s}\big )-\rV_{\rT}\big (\infty;\xh_{s}\big )&\hspace*{-.2cm}=
	\dfrac{\big (\rC_{\rT}(x_{s})+\rG_{\rE}(\rV_{\rE}(\infty))-\rC_{\rT}(\xh_{s})-\rG_{\rE}(\rVh_{\rE}(\infty))\big )}{1+
		\rG_{\rI}'(\xi_{\rI})\rG_{\rT}'(\xi_{\rT})}\\
	&
	-\dfrac{\rG_{\rE}(\rVh_{\rE}(\infty))\big )-\rG_{\rI}'(\xi_{\rI})\big (\rC_{\rI}(x_{s})-\rC_{\rI}(\xh_{s})\big )}{1+
		\rG_{\rI}'(\xi_{\rI})\rG_{\rT}'(\xi_{\rT})},
\end{array}
$$ 
for some intermediate values $\xi_{\rI}$ y $\xi_{\rT}$. Since $\rC_{\rI}$ does not depend on  $x_{s}$ (see \eqref{eq:constantsterms}) one has
$$
\rV_{\rT}\big (\infty;x_{s}\big )-\rV_{\rT}\big (\infty;\xh_{s}\big )=
\dfrac{g_{s\rT}(x_{s})+\rG_{\rE}\big (\rV_{\rE 0}+g_{s\rE}(x_{s})\big )-g_{s\rT}(\xh_{s})-\rG_{\rE}\big (\rV_{\rE 0}+g_{s\rE}(\xh_{s})\big )}{1+\rG_{\rT}'(\xi_{\rT})\rG_{\rI}'(\xi_{\rI})}.
$$ 
Hence, the monotonicity of  $\rC_{\rT}(x_{s})$ is transferred to $\rV_{\rT}(\infty;x_{s})$.
\par
\noindent
b) Analogously,  \eqref{eq:variacionVt} leads to
$$
\rV_{\rT}\big (\infty;\alpha_{1}\big )-\rV_{\rT}\big (\infty;\alpha_{2}\big )=
\dfrac{\big (\rC_{\rT}(\alpha_{1})-\rC_{\rT}(\alpha_{2})\big )-\rG_{\rI}'(\alpha_{\rI})\big (\rC_{\rI}(\alpha_{1})-\rC_{\rI}(\alpha_{2})\big )}{1+
	\rG_{\rI}'(\alpha_{\rI})\rG_{\rT}'(\alpha_{\rT})},
$$ 
for some intermediate values $\alpha_{\rI}$ y $\alpha_{\rT}$, where the potentials $\rV_{\rI}\big (\infty;\alpha_{1}\big ), \rV_{\rI}\big (\infty;\alpha_{2}\big )$ and
$\rV_{\rT}\big (\infty;\alpha_{1}\big )$, $\rV_{\rT}\big (\infty;\alpha_{2}\big )$, as well as   $\rC_{\rI}(\alpha_{1}), \rC_{\rI}(\alpha_{2})$ and
$\rC_{\rT}(\alpha_{1}), \rC_{\rT}(\alpha_{2})$, are given in the corresponding system~\eqref{eq:sistemaestacionarionotacional}. Thus
$$
\rV_{\rT}\big (\infty;\alpha_{1}\big )-\rV_{\rT}\big (\infty;\alpha_{2}\big )=-
\dfrac{\rG_{\rI}'(\alpha_{\rI})g_{c\rI}\big [\psi (x_{l})\big ]}{1+
	\rG_{\rI}'(\alpha_{\rI})\rG_{\rT}'(\alpha_{\rT})}(\alpha_{1}-\alpha_{2}),
$$ 
because $\rC_{\rT}$ is constant with respect to $\alpha_{c\rI}$ (see again \eqref{eq:constantsterms}). Now, the monotonicity of the right-hand side with respect to $\alpha_{c\rI}$ is transferred to $\rV_{\rT}(\infty;\alpha_{c\rI})$ in the opposite direction. $\fin$
\par
\begin{figure}[ht]
\centering
\includegraphics[width=10.5cm]{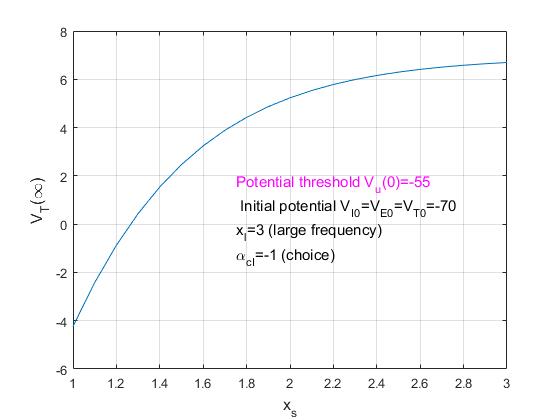} 
\caption{ Monotone dependence of $\rV_{\rT}$ with respect to $x_{s}$}
\label{fig:VariacionVtfrecuenciacorta}
\end{figure}
\begin{figure}[ht]
\centering
\includegraphics[width=10.5cm]{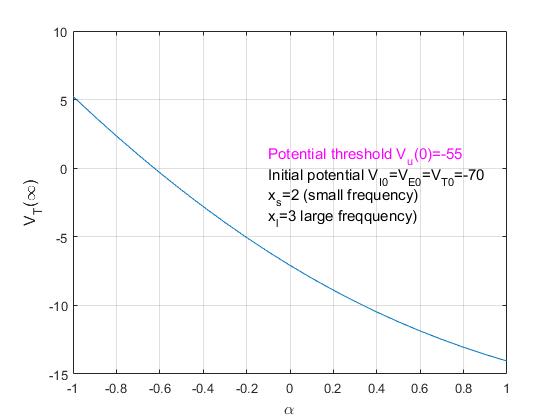} 
\caption{Monotone dependence of $\rV_{\rT}$ with respect to $\alpha_{c\rI}$}
\label{fig:VariacionVtmodo}
\end{figure}

\bigskip
As it was pointed out in Section \ref{sec:model}, Figures \ref{fig:VariacionVtfrecuenciacorta} and \ref{fig:VariacionVtmodo} illustrate the above monotone dependence 
\begin{rem}\rm Item a) of Proposition \ref{prop:variationVtparameters} corresponds to comment (i) made in \cite{BS} and quoted at the beginning of Section 2. Item b) of Proposition \ref{prop:variationVtparameters} corresponds to comments (iv) and (v) made in~\cite{BS} and quoted at the beginning of Section 2. $\fin$
\end{rem}
\begin{rem}\rm
The biological interest of the above result is quite evident:
\begin{enumerate}
\item Part a) of Proposition \ref{prop:variationVtparameters} extends Lemma 5 of \cite{BS}. Inequality \eqref{eq:Vtxs} indicates that, under constant pain, when the stimulation of the short fibers is slightly increased without any other change, then, after a while, the pain felt will again be constant but of greater intensity.
\item Part b) of Proposition \ref{prop:variationVtparameters} extends Lemma 8 of \cite{BS}. Inequality \eqref{eq:Vtaci} indicates that activating the cognitive control mechanism either reduces or increases the equilibrium state of T-cell potentials, depending on whether it has an inhibitory or excitatory effect. Again, assuming that the solution of the system is the equilibrium solution, this would correspond to observations (iv) and (v) of \cite{BS}. $\fin$
\end{enumerate}
\label{obser:variacionVtparametros}
\end{rem}
\begin{rem}\rm Given $\rV_{\rE 0}$, we can study the relative potentials  $\rV_{\rI}\big (\infty;x_{l}\big ),
\rV_{\rI}\big (\infty;\xh_{l}\big )$ and
$\rV_{\rT}\big (\infty;x_{l}\big )$, $\rV_{\rT}\big (\infty;\xh_{l}\big ) $, as well as the expressions $\rC_{\rI}(x_{l}), \rC_{\rI}(\xh_{l})$ and $\rC_{\rT}(x_{l})$, $\rC_{\rT}(\xh_{l})$ given in  the corresponding system 
\eqref{eq:sistemaestacionarionotacional}. So, for $\alpha_{c\rI}=0$, given $x_{s}$, formula \eqref{eq:variacionVt} leads to
\begin{equation}
\rV_{\rT}\big (\infty;x_{l}\big )-\rV_{\rT}\big (\infty;\xh_{l}\big )=
\dfrac{\big (g_{l\rT}-\rG_{\rI}'(\xi_{\rI})g_{l\rI}\big )(x_{l})-\big (g_{l\rT}-\rG_{\rI}'(\xi_{\rI})g_{l\rI}\big )(\xh_{l})\big )}{1+ \rG_{\rI}'(\xi_{l})\rG_{\rT}'(\xi_{l})}.
\label{eq:Vtxl}
\end{equation}
For some intermediate values $\xi_{\rI}$ and $\xi_{\rT}$, we can use the monotonicity of the functions $g_{l\rT}-\rG_{\rI}'(\xi_{\rI}) g_{l\rI}$ with respect to $x_{l}$ in order to analyze the dependence of $\rV_{\rT}\big (\infty;x_{l}\big )$ on $x_{l}$. Inequality \eqref{eq:Vtxl} extends Lemma 6 of \cite{BS}. Note also that inequality \eqref{eq:Vtxl} indicates that if no type of cognitive control is exercised, constant pain is felt, and the stimulation of the long fibers is slightly increased without other changes; then, after a short time, the pain felt will also be constant, but of greater or lesser intensity. This will depend on the details of the model and the level of stimulation of the short and long fibers considered. Therefore, the Gate Control Theory of Pain can explain different consequences following an increase in the excitation of long fibers. In particular, the possibility that the pain increases temporarily until it stabilizes at a lower intensity would explain observation (ii) of \cite{BS} quoted in Section 2. It is also interesting to note that Nathan and Rudge \cite{NR} found that stimulating the long fibers did not always reduce the pain caused by the short fibers: they used this fact as an argument against the Gate Control Theory of Pain. However, as we have shown, this theory is perfectly capable of explaining such observations. $\fin$
\label{obser:Vtxl}
\end{rem}
\begin{rem}\rm There are many possible variants. For instance, we can maintain constant the parameters $\rV_{\rE 0}, x_{s}$, $x_{l}$ and $\alpha_{c\rI}$ in order to study the dependence with respect to the function $f_{\rE}$. Under such conditions, we consider the potentials  $\rV_{\rI}\big (\infty;f_{\rE}\big )$, $\rV_{\rI}\big (\infty;\widehat{f}_{\rE}\big )$ and $\rV_{\rT}\big (\infty;f_{\rE}\big ), \rV_{\rT}\big (\infty;\widehat{f}_{\rE}\big )$, as well as the expressions $\rC_{\rI}(f_{\rE}), \rC_{\rI}(\widehat{f}_{\rE})$ and $\rC_{\rT}(f_{\rE}), \rC_{\rT}(\widehat{f_{\rE}})$ given in the corresponding system ~\eqref{eq:sistemaestacionarionotacional}. Since 
$\rC_{\rI}-\rC_{\rI}=0$ and  $\rC_{\rT}-\rC_{\rT}=0,$
formula \eqref{eq:variacionVt} leads to
\begin{equation}
\rV_{\rT}\big (\infty;f_{\rE}\big )-\rV_{\rT}\big (\infty;\widehat{f}_{\rE}\big )=
\dfrac{g_{\rE\rT}\big [f_\rE(\rV_\rE(\infty)) \big ]-g_{\rE\rT}\big [\widehat{f}_\rE(\rV_\rE(\infty)) \big ]}{1+\rG_{\rI}'(\xi_{\rI})\rG_{\rT}'(\xi_{\rT})},
\label{eq:Vtfe}
\end{equation}
for some intermediate values $\xi_{\rI}$ y $\xi_{\rT}$.  Then
$$
f_\rE(\rV_\rE(\infty))\le \widehat{f}_\rE(\rV_\rE(\infty)) \quad \Rightarrow \quad \rV_{\rT}\big (\infty;,f_{\rE}\big )\le \rV_{\rT}\big (\infty,\widehat{f}_{\rE}\big )
$$
follows. $\fin$
\label{obser:Vtfe}
\end{rem}

Our last result on the dependence of the potential  $\rV_{\rT}$ concerning the function $\varphi$. We emphasize \eqref{eq:equlibriumexistence} by
\begin{equation}
\rV_{\rT}+\rG_{\rI}\bigg (\overbrace{g_{d\rI}\big (\varphi\big [f_\rT(\rV_\rT)\big ]\big )}^{\rG_{\rT}^{\varphi}(\rV_{\rT})}\bigg )+\rC_{\rI}\big )=
\overbrace{\rC_{\rT}+\rG_{\rE}\big (\rV_{\rE}(\infty)\big )}^{\hbox{\tiny (independent of the potentials $\rV_{\rI}$ and $\rV_{\rT}$)}}
\label{eq:equlibriumexistence1}
\end{equation}
(see \eqref{eq:Gfunctions}). We recall the action of $\varphi$ in the modeling given in \eqref{eq:frecuencies}.

\begin{prop}[Dependence of $\rV_{\rT}(\infty)$ with respect to $\varphi$] Assume \eqref{eq:constant in time}. Given $\rV_{\rE 0}$,  $x_{l}$,  $x_{l}$ and $\alpha _{c\rI}$,  the stationary potential $\rV_{\rT}(\infty;\varphi)$ depends monotonically but in the opposite sense to $\varphi.\fin$ 
\label{prop:variationVtphi}
\end{prop}
\noindent 
\proof Given the functions $\varphi$ and $\widehat{\varphi}$, we consider the associated potentials $\rV_{\rT}\big (\infty;\varphi\big )$ and $\rV_{\rT}\big (\infty;\widehat{\varphi}\big )$. Then 
$$
\rV_{\rT}\big (\infty;\varphi\big )+\rG_{\rI}\big (\rG_{\rT}^{\varphi}\big (\rV_{\rT}\big (\infty;\varphi\big )\big )+\rC_{\rI}\big )=\rV_{\rT}\big (\infty;\widehat{\varphi}\big )+\rG_{\rI}\big (\rG_{\rT}^{\widehat{\varphi}}\big (\rV_{\rT}\big (\infty;\widehat{\varphi}\big )\big )+\rC_{\rI}\big ).
$$
holds (see  \eqref{eq:equlibriumexistence1}). By construction,  we have
$$
\rG_{\rT}^{\varphi}(\rV)>\rG_{\rT}^{\widehat{\varphi}}(\rV),
$$
whenever $\varphi>\widehat{\varphi}$. So that the monotonicity of the function 
$$
\rV\mapsto \rV+\rG_{\rI}\big (\rG_{\rT}^{\varphi}\big (\rV\big )\big )
$$
concludes the result $.\fin$

\begin{figure}[ht]
\centering
\includegraphics[width=10.5cm]{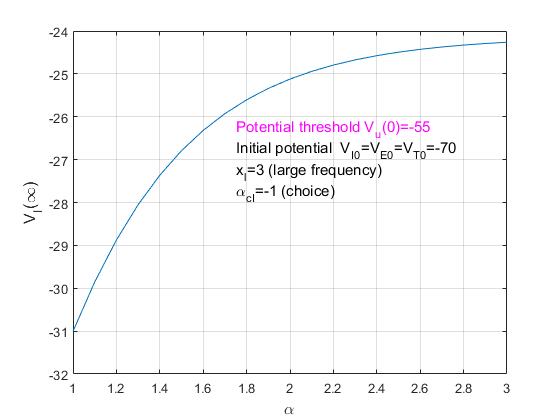} \\ [-.35cm]
\caption{Dependence of $\rV_{\rI}$ with respect to $x_{s}$}
\label{fig:VariacionVifrecuenciacorta}.
\end{figure}

\begin{figure}[ht]
\centering
\includegraphics[width=9cm]{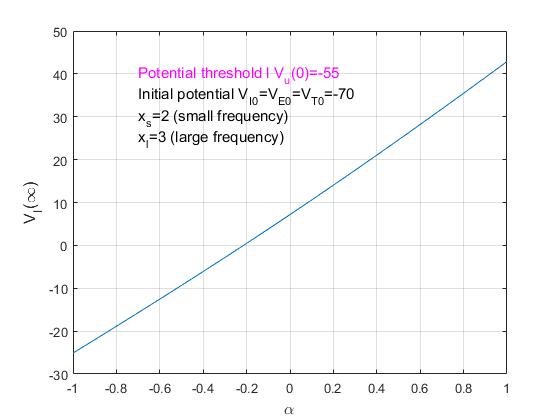} \\ [-.35cm]
\caption{Dependence of $\rV_{\rI}$ with respect to   $\alpha_{c\rI}$}.
\label{fig:VariacionVImodo}
\end{figure}

\begin{rem}\rm 
Proposition \ref{prop:variationVtphi} extends Lemma 7 of \cite{BS}. The consequence of this result is that we can reduce the equilibrium point value of the T-cell potential by increasing the signal from the mid-brain (the inhibitory step-down control). Assuming that the system is in a state of equilibrium, this shows that the pain perception is more cushioned (observation (iii) of \cite{BS} quoted at the beginning of Section 2). $\fin$
\label{obser:variacionVtphi}
\end{rem}
We can get an analogous version of Proposition \ref{prop:variationVtparameters} but now for the potential $\rV_{\rI}(\infty)$: it follows again from \eqref{eq:variacionVi}.
\begin{prop}[Dependence  $\rV_{\rI}(\infty)$ with respect to the data] 
\
\par
\noindent 
Assume \eqref{eq:constant in time}.
\par
\noindent 
a) {\sc Dependence on $\alpha_{c\rI}$}.
Given $\rV_{\rE 0}$,  $x_{s}$ and $x_{l}$,  the stationary potential $\rV_{\rI}(\infty;\alpha_{c\rI})$ depends increasingly on $\alpha_{c\rI}$. Moreover, we have the inequality
\begin{equation}
\big [\rV_{\rI}\big (\infty;\alpha_{1}\big )-\rV_{\rI}\big (\infty;\alpha_{2}\big )\big ]_{+}\le 
g_{c\rI}\big [\psi (x_{l})\big ]\big [\alpha_{1}-\alpha_{2}\big ]_{+}.
\label{eq:Viaci}
\end{equation}
b) {\sc Dependence on $x_{s}$}. Given $x_{l}$ and $\alpha_{ci}$ the stationary potential $\rV_{\rI}(\infty;x_{s})$ depends monotonically on $x_{s}$ but reversing the monotonicity of the expression $g_{s\rT}(x_{s})+g_{\rE\rT}\big [f_\rE(\rV_{\rE 0} + g_{s\rE}(x_s)) \big ].\fin$
\label{prop:Vi}
\end{prop} 
\proof Given $\alpha_{1}$ and $\alpha_{2}$,  the argument used in the proof of Proposition \ref{prop:variationVtparameters},  from \eqref{eq:variacionVi}, we obtain
$$
\rV_{\rI}\big (\infty;\alpha_{1}\big )-\rV_{\rI}\big (\infty;\alpha_{2}\big )=
\dfrac{\big (\rC_{\rI}(\alpha_{1})-\rC_{\rI}(\alpha_{2})\big )-\rG_{\rT}'(\alpha_{\rT})\big (\rC_{\rT}(\alpha_{1})-\rC_{\rT}(\alpha_{2})\big )}{1+\rG_{\rI}'(\alpha_{\rI})\rG_{\rT}'(\alpha_{\rT})},
$$ 
for some intermediate values $\alpha_{\rT}$ and $\alpha_{\rI}$. Since  $\rC_{\rT}(\alpha_{c\rI})$ does not depend $\alpha_{c\rI}$ (see \eqref{eq:constantsterms}) we get
$$
\rV_{\rI}\big (\infty;\alpha_{1}\big )-\rV_{\rI}\big (\infty;\alpha_{\rI}\big )=
\dfrac{g_{c\rI}\big [\psi (x_{l})\big ]}{1+\rG_{\rT}'(\alpha_{\rT})\rG_{\rI}'(\alpha_{\rI})}(\alpha_{1}-\alpha_{2}).
$$ 
Then the monotonicity of $\rC_{\rI}(\alpha_{c\rI})$ with respect to $\alpha_{c\rI}$ is transferred to the function $\rV_{\rI}(\infty;\alpha_{c\rI})$.
\par
\noindent
b) The dependence on $x_{s}$ results from the identity
$$
\begin{array}{ll}
	\rV_{\rI}\big (\infty;x_{s}\big )-\rV_{\rI}\big (\infty;\xh_{s}\big )&\hspace*{-.2cm}=
	\dfrac{\big (\rC_{\rI}(x_{s})-\rC_{\rI}(\xh_{s})\big )}{1+\rG_{\rI}'(\xi_{\rI})\rG_{\rT}'(\xi_{\rT})}\\ [.3cm]
	&\hspace*{-3cm}
	-\dfrac{\rG_{\rT}'(\xi_{\rT})}{1+\rG_{\rI}'(\xi_{\rI})}\rG_{\rT}'(\xi_{\rT})\big [g_{s\rT}(x_{s})+\rG_{\rE}(\rV_{\rE 0}+g_{s\rE}(x_s))-g_{s\rT}(\xh_{s})-\rG_{\rE}(\rV_{\rE 0}+g_{s\rE}(\xh_s))\big ],
\end{array}
$$ 
for some intermediate values $\xi_{\rT}$ and $\xi_{\rI}$. Thus
$$
\rV_{\rI}\big (\infty;x_{s}\big )-\rV_{\rI}\big (\infty;\xh_{s}\big )=-
\dfrac{\rG_{\rT}'(\xi_{\rT})\big [g_{s\rT}(x_{s})+\rG_{\rE}(\rV_{\rE 0}+g_{s\rE}(x_s))-g_{s\rT}(\xh_{s})-\rG_{\rE}(\rV_{\rE 0}+g_{s\rE}(\xh_s))\big ]}{1+\rG_{\rI}'(\xi_{\rI})\rG_{\rT}'(\xi_{\rT})},
$$
since $\rC_{\rI}(x_{s})$ does not depend on $x_{s}$ (see again \eqref{eq:constantsterms}). Therefore the monotonicity of  $x_{s}\mapsto g_{s\rT}(x_{s})+g_{\rE\rT}\big [f_\rE(\rV_{\rE 0} + g_{s\rE}(x_s)) \big ]$ is transferred to the one of the function $\rV_{\rI}(\infty;x_{s})$ but in the opposite sense$.\fin$ 
\begin{rem}\rm
Note that the function $x_{s} \mapsto g_{s\rT}(x_{s}) + g_{\rE\rT}\left[ f_\rE(\rV_{\rE 0} + g_{s\rE}(x_s)) \right]$ may be strictly decreasing for some negative values of the stationary potential $\rV_{\rE 0}$, even though the functions $g_{s\rT}(x)$, $g_{\rE\rT}(x)$, and $f_\rE(x)$ are strictly increasing. This occurs in Figure \ref{fig:VariacionVifrecuenciacorta} that illustrates part b) of Proposition \ref{prop:Vi}. On the other hand, 
Figure \ref{fig:VariacionVImodo} illustrates part a) of Proposition  \ref{prop:Vi}.$\fin$ 
\label{obser:Vi}
\end{rem} 

\begin{rem}\rm It is possible to get some monotone dependence results similar to the ones quoted in Remarks \ref{obser:Vtxl} and \ref{obser:Vtfe}.$\fin $
\label{obser:Vi1}
\end{rem}

\begin{rem}\rm
The above results remain valid, with slight changes, when we replace a constant small frequency $x_{s}$ by a bounded time-dependent frequency, $x_{s}^{v}(t),$ verifying
\begin{equation}
\lim_{t\rightarrow \infty}x_{s}^{v}(t)=x_{s}.
\label{eq:xsvariableinfinito}
\end{equation}
In that case, the uncoupled equation
$$
\tau_\rE \derp{\rV}_\rE(t) =-\big (\rV_\rE(t)-\rV_{\rE 0}\big ) + g_{s\rE}\big (x_s(t)\big ) 
$$
leads to
\begin{equation}
\rV_{\rE}(t)=\big (\rV_{\rE}(0)-\rV_{\rE 0}\big )e^{-\frac{t}{\tau_{\rE}}}+ \rV_{e0}+\dfrac{1}{\tau_{\rE}}\int_{0}^{t}g_{s\rE}(x_{s}^{v}(r))e^{-\frac{t-r}{\tau_{\rE}}}dr,\quad t\ge 0. 
\label{eq:Veinfinitovariable}
\end{equation}
Since $g_{s\rE}(x_{s}^{v}(t))\nearrow g_{s\rE}(x_{s})$, as $t\nearrow \infty$, we deduce
$$
\dfrac{1}{\tau_{\rE}}\int_{0}^{t}g_{s\rE}(x_{s}^{v}(r))e^{-\frac{t-r}{\tau_{\rE}}}dr-
g_{s\rE}(x_{s})\big (1-e^{-\frac{t}{\tau_{\rE}}}\big )=
\dfrac{1}{\tau_{\rE}}\int_{0}^{t}\big (g_{s\rE}(x_{s}^{v}(r))e^{-\frac{t-r}{\tau_{\rE}}}- g_{s\rE}(x_{s})\big )e^{-\frac{t-r}{\tau_{\rE}}}dr
$$
and then
$$
\left |\dfrac{1}{\tau_{\rE}}\int_{0}^{t}g_{s\rE}(x_{s}^{v}(r))e^{-\frac{t-r}{\tau_{\rE}}}dr-
g_{s\rE}(x_{s})\big (1-e^{-\frac{t}{\tau_{\rE}}}\big )\right |\le \varepsilon \big (1-e^{-\frac{t}{\tau_{\rE}}}\big )<\varepsilon,
$$
assumed $|g_{s\rE}(x_{s}^{v}(r))- g_{s\rE}(x_{s})|\le \varepsilon$. Then we get that
$\rV_{\rE}(\infty)\doteq  \rV_{\rE 0} + g_{s\rE}(x_s)$ which coincides with the value obtained in the case of constant frequencies  (see \eqref{eq:Veinfinito}). $\fin$
\end{rem}
\par
\begin{rem}\rm
As a consequence, all the above results concerning the stationary states remain valid for bounded time-dependent frequencies satisfying \eqref{eq:xsvariableinfinito}. $\fin$
\end{rem}

\begin{figure}[ht]
\centering
\includegraphics[width=10.5cm]{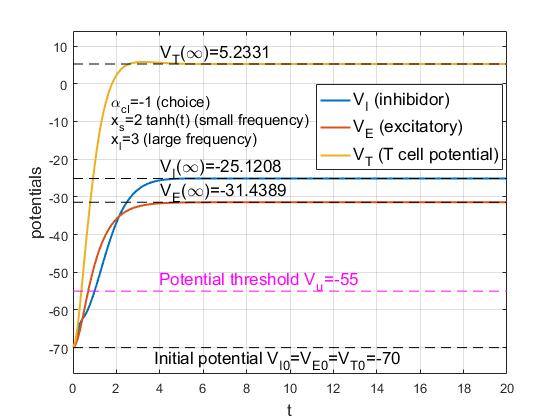} \\ [-.35cm]
\caption{Representation of the potentials for time-depending frequencies}
\label{fig:Graficassolucionvariable}
\end{figure}

\begin{figure}[ht]
\centering
\includegraphics[width=10.5cm]{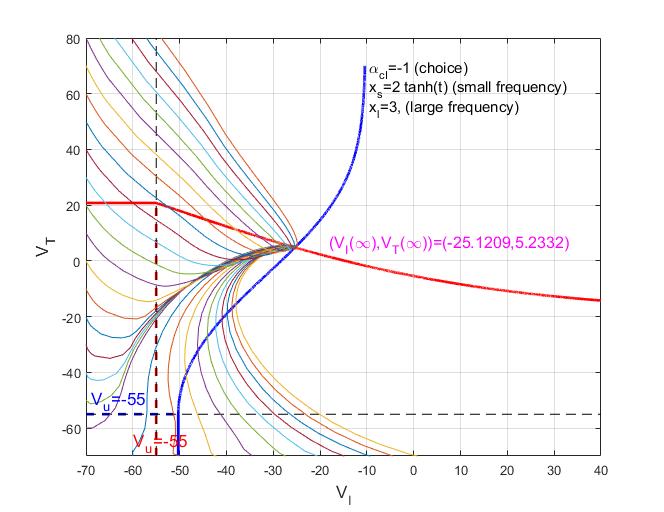} \\ [-.35cm]
\caption{Phase plane for some variable frequencies}
\label{fig:planofasesvariable}
\end{figure}

\begin{figure}
\centering
\includegraphics[width=11cm]{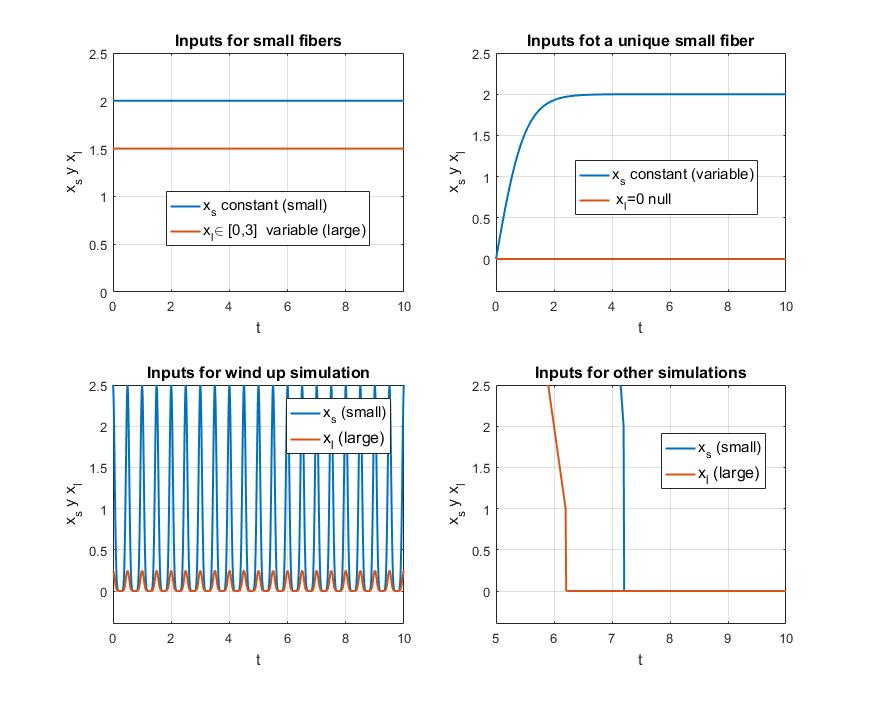} \\ [-.35cm]
\caption{Other input types of frequencies}
\label{fig:inputs.jpg}
\end{figure}

\begin{rem}\rm
Following \cite{BCS}, we can consider other types of frequency inputs (see Figure \ref{fig:inputs.jpg}). The mathematical analysis raises the intriguing possibility of oscillatory solutions to the equations, as suggested by the observation behind Lemma 4 of \cite{BS}. If such a solution exists, the $\rV_\rT$ potential of the T-cells oscillates, causing the pain to increase and decrease rhythmically. In this case, the model prediction would be that the transition from constant to rhythmic pain can only occur due to a drastic variation in firing rates in the small or large fibers, assuming no change in the descending controls. $\fin$
\end{rem}

\section{Optimal Control Problem}
\label{sec:controlproblem}
In this section, we consider the optimal control problem associated with the equation \eqref{eq:sistemadinamico} for the potentials $\bV(t) = \big (\rV_{\rE}(t), \rV_{\rI}(t), \rV_{\rT}(t)\big )^{\tt t}$ when they are controlled by the short frequency $x_{s}$. Let us emphasize the system in the following way:
\begin{equation}
\left \{
\begin{array}{l}
\derp{\bV} (t)=\bff \big (\bV(t),x_{s}(t)\big ),\quad 0\le t,\\ [.175cm]
\bV(0)=\bV_{0},
\end{array}
\right .
\label{eq:controlledsystem}
\end{equation}
where $ \bV_{0}=(\rV_{\rE 0},\rV_{\rI 0},\rV_{\rT 0})^{\tt t}$ and
\begin{equation}
\bff \big (\bV,x_{s}\big )=
\left (
\begin{array}{l}
\big (-\rV_\rE+\rV_{\rE 0} + g_{s\rE}(x_s)\big )\tau_\rE ^{-1} \\ [.1cm]
\big (-\rV_\rI+\rG_{\rT}(\rV_{\rT})+\rC_{\rI}\big )\tau_\rI^{-1} \\ [.1cm]
\big (-\rV_\rT-\rG_{\rI}(\rV_{\rI})+\rG_{\rE}\big (\rV_{\rE}\big )+\underbrace{\rV_{\rT 0} + g_{s\rT}(x_s) + g_{l\rT}(x_l)}_{\rC_{\rT}}\big )\tau_\rT^{-1}
\end{array}
\right )
\end{equation}
(see  \eqref{eq:field}, \eqref{eq:Gfunctions} and  \eqref{eq:constantsterms}). 
The formulation as an optimal control problem comes from the minimization of the functional in which the pain at the final time $t_{f} > 0$ is assumed to be an increasing function of the $\rV_{\rT}(t_{f})$ potential
\begin{equation}
\rW(\bV_{0})\doteq \min_{x_{s}\in \cU}\J_{\bV_{0}} (x_{s}),\quad \J_{\bV_{0}} (x_{s})\doteq \int^{t_{f}}_{0}\left (
\dfrac{1}{2}|\rV_{\rT}(\sigma)|^{2}- \rQ g_{s\rT}\big (x_{s}(\sigma)\right )d\sigma+\rS\big (\rV_{\rT}(t_{f})\big)
\label{eq:minimization1}
\end{equation}
(see \eqref{eq:minimization}). Here, the set of controls is given by $\cU=\rL^{\infty}\big (\R_{+}:[\underline{x_{s}},\overline{x_{s}}]\big )$, $\rQ$ is a given positive constant and $\rS$ a $\cC^{1}$ function satisfying {\em $0<\rS_{-}'\le \rS'\big (\rV_{\rT}\big )\le \rS'_{+}$}. We emphasize that the control values belong to the closed interval $[\underline{x_{s}},\overline{x_{s}}]\subset \R$.
\par
As already mentioned, from \cite[Theorem 3.2.1]{BP}, for each $x_{s}\in \rL^{1}(0,t_{f};\R)$, the system \eqref{eq:controlledsystem} admits a unique solution $\bV(t)$, absolutely continuous, defined in the interval $[0,t_{f}]$. Thus, $\bV \in \cC\big([0,t_{f}];\R^{3}\big) \cap \rW^{1,1}(0,t_{f};\R^{3})$. Moreover, as proved in \cite{BP}, $\bV(t)$ depends continuously on $x_{s}(\cdot)$, so that
$
\bV \in\cC \big ( \rL^{1}(0,t_{f}:\R):\cC\big ([0,t_{f}]:\R^{3}\big )\big ).
$ 
\par
Following \cite{FR} we consider the states, $\bV$, the  co-states , $\bP$,  as well as the  optimality Hamiltonian
$$
\rH(\bV,\bP,a)=\pe{\bff (\bV,a)}{\bP}-\dfrac{1}{2}|\rV_{\rT}|^{2}+\rQ g_{s\rT}(a),\quad (\bV,\bP)\in\R^{3}\times \R^{3},~a\in [\underline{x},\overline{x}].
$$
The optimality Hamiltonian system is given by
\begin{equation}
\left \{
\begin{array}{ll}
\derp{\bV^{*}}(t)=\rD_{\bP}\rH \big (\bV^{*}(t),\bP^{*}(t),x_{s}^{*}(t)\big ),&\quad \bV^{*}(0)=\bV_{0}\in\R^{3},\\ [.125cm]
\derp{\bP^{*}}(t)=-\rD_{\bV}\rH \big (\bV^{*}(t),\bP^{*}(t),x_{s}^{*}(t)\big ),&\quad \bP^{*}(t_{f})=\big (0,0,-\rS'\big (\rV_{\rT}^{*}(t_{f}\big )\big )^{\tt t}\in\R^{3}.
\end{array}
\right .
\label{eq:ecHamilton}
\end{equation}
\vspace*{.1cm}
By \cite[Theorem 6.1.1]{BP}, we obtain the existence and uniqueness of the solution $\big(\bV^{*}(t), \bP^{*}(t)\big )^{\tt t}$ of the Hamilton system \eqref{eq:ecHamilton}, as well as the existence of an optimal control $x_{s}^{*} \in \rL^{\infty}(0,t_{f};[\underline{x_{s}}, \overline{x_{s}}])$, which may be discontinuous (see also Theorem \ref{theo:existenceoptimalcontrol} below).
\par
Obviously, the first equation in \eqref{eq:ecHamilton} coincides with \eqref{eq:controlledsystem} at $x_{s}=x_{s}^{*}$. On the other hand,  the adjoint equation in \eqref{eq:ecHamilton} is the  linear system on $\bP^{*}(t)=(\rP^{*}_{\rE}(t),\rP^{*}_{\rI}(t),\rP^{*}_{\rT}(t)\big )^{\tt t}$
\begin{equation}
\derp{\bP^{*}}(t)=\cA(t) \bP^{*}(t) +\bF (t),
\label{eq:moments}
\end{equation}
with
$$
\hspace*{-.2cm}
\cA(t)=\left (
\begin{array}{ccc}
 \tau_\rE ^{-1}	& &\\ [.15cm]
 &  \tau_{\rI}^{-1} & -\tau_{\rI}^{-1}\rG_{\rT}'\big (\rV_{\rT}^{*}(t)\big )\\ [.1cm]
\tau_{\rT}^{-1}\rG_{\rI}'\big (\rV_{\rI}^{*}(t)\big ) &-\tau_{\rT}^{-1}\rG_{\rE}'\big (\rV_{\rE}^{*}(t)\big )&\tau_{\rT}^{-1}
\end{array}
\right )\quad \hbox{and}\quad 
\bF(t)=-
\left (
\begin{array}{c}
0\\ 
0\\
\rV^{*}_{\rT}(t)
\end{array}
\right )
$$
(see \eqref{eq:ecHamilton}). Whenever $\bP_{1}^{*}(t),\bP_{2}^{*}(t) $ and $\bP_{3}^{*}(t)$ are   linearly independent solutions of the homogeneous problem
\begin{equation}
\derp{\bP^{*}}(t)=\cA(t) \bP^{*}(t),\quad t <t_{f},
\label{eq:momentsunforced}
\end{equation}
we may construct the fundamental matrix 
$$
\Phi(t-t_{f})\doteq
\left (
\begin{array}{ccc}
\rP_{\rE 1}^{*}(t) & \rP_{\rE 2}^{*}(t) & \rP_{\rE 3}^{*}(t)\\ [.15cm]
\rP_{\rI 1}^{*}(t) & \rP_{\rI 2}^{*}(t) & \rP_{\rI 3}^{*}(t)\\ [.15cm]
\rP_{\rT 1}^{*}(t) & \rP_{\rT 2}^{*}(t) & \rP_{\rT 3}^{*}(t)
\end{array}
\right ).
$$
Then from the constant variation formula the solution of \eqref{eq:moments} is given by 
\begin{equation}
\bP^{*}(t)=\Phi \big (-(t_{f}-t)\big )\left (\Phi ^{-1}(0)\bP^{*}(t_{f})+\int^{t_{f}}_{t}\Phi ^{-1}\big (-(t_{f}-s)\big ) \bF(s)ds\right ),\quad t <t_{f}.
\label{eq:variacionconstanteP}
\end{equation}
\par
We note that the equation for the optimal co-state $\rP_{\rE}^{*}(t)$ satisfies 
$$
\derp{\rP_{\rE}^{*}}(t)=\tau_\rE ^{-1}\rP_{\rE}^{*}(t),\quad \rP_{\rE}^{*}(t_{f})=0,
$$
whence
\begin{equation}
\rP_{\rE}^{*}(t)\equiv 0,\quad t \le t_{f}.	\label{eq:momentE}
\end{equation}
Then \eqref{eq:moments} becomes the  reduced linear system for $\widehat{\bP}^{*}(t)=(\rP^{*}_{\rI}(t),\rP^{*}_{\rT}(t)\big )^{\tt t}$ 
\begin{equation}
\derp{\widehat{\bP}^{*}}(t)=\widehat{\cA(t)} \widehat{\bP}^{*}(t) +\widehat{\bF} (t),
	\label{eq:hatmoments}
\end{equation}
with
$$
\widehat{\cA}(t)=\left (
\begin{array}{cc}
\tau_{\rI}^{-1}	 & -\tau_{\rI}^{-1}\rG_{\rT}'\big (\rV_{\rT}^{*}(t)\big )\\ [.1cm]
	-\tau_{\rT}^{-1}\rG_{\rE}'\big (\rV_{\rE}^{*}(t)\big ) &\tau_{\rT}^{-1}
\end{array}
\right )\quad \hbox{and}\quad 
\widehat{\bF}(t)=-
\left (
\begin{array}{c}
	0\\
	\rV^{*}_{\rT}(t).
\end{array}
\right ).
$$

\begin{theo}[Existence of an Optimal Control] Let us assume that 
\begin{equation}
	\hbox{the functions   $\rG_{\rT}(\rV)$ and $\rG_{\rE}(\rV)$ are convex and the function $\rG_{\rI}(\rV)$  is concave}
	\label{eq:convexity}
\end{equation}
{\rm (}see  \eqref{eq:Gfunctions}{\rm )}. Then there exists an optimal control $x_{s}^{*}\in \rL^{\infty}(\underline{x},\overline{x})$.
\label{theo:existenceoptimalcontrol} 
\end{theo}
\proof The result follows from the proof of  Theorem 5.2.2 of \cite{BP}. Essentially, it is based on the  convexity of the set
$$
\cF(\bV)=\left \{(\by,y_{0})\in\R^{4}:~y_{0}\ge \dfrac{1}{2}|\rV_{\rT}|^{2}+\rQ g_{s\rT}(x_{s}),~\by(t)=\bff \big (\bV(t),x_{s}\big ),~x_{s}\in[\underline{x},\overline{x}] \right \},
\quad \bV\in\R^{3}.
$$
$\fin$
\begin{rem}\rm We emphasize that the above proof does not require the convexity of the functional $\J_{\bV_{0}}(x_{s})$ with respect to the control $x_{s}(\cdot)$. This remark was already pointed out in \cite{Tr}. $\fin$

\end{rem}
\par
The property $\rP_{\rE}^{*}(t)\equiv 0$ is very useful. Indeed, the Hamiltonian applied to the  optimal states becomes
$$
\begin{array}{ll}
\rH\big (\bV^{*},\bP^{*},a\big )&\hspace*{-.2cm}=\left (-\rV_\rI^{*}+\rG_{\rT}(\rV_{\rT}^{*})+\rC_{\rI}\right )\tau_\rI^{-1}\rP^{*}_{\rI}\\ 
&+\left (-(\rV_\rT^{*}-\rV_{\rT 0}) + g_{s\rT}(a) + g_{l\rT}(x_l) +\rG_{\rE}\big (\rV_{\rE}^{*}\big ) -\rG_\rI\big (\rV_\rI^{*}\big )\right )\tau_\rT^{-1}\rP^{*}_{\rT} \\
&
-\dfrac{1}{2}(\rV_\rT^{*}(\sigma))^{2}+\rQ g_{s\rT}(a).
\end{array}
$$
Therefore the Pontryagin Maximum Principle (see \cite{FR}) implies that
$$
\rH(\bV^{*}(t),\bP^{*}(t),x_{s}^{*}(t))=\max_{a	\in[\underline{x_{s}},\overline{x_{s}}]}\rH(\bV^{*}(t),\bP^{*}(t),a).
$$
A reduced version of it leads to
$$
\big (\tau_\rT^{-1}\rP^{*}_{\rT}(t)-\rQ \big )g_{s\rT}(x_{s}^{*}(t))=\max_{a	\in[\underline{x_{s}},\overline{x_{s}}]}\left [\big (\tau_\rT^{-1}\rP^{*}_{\rT}(t)+\rQ \big )g_{s\rT}(a)\right ].
$$
This shows that the optimal control $x_{s}^{*}(t)$ is determined by the maximum values of the bounded real function.
\begin{equation}
[\underline{x_{s}},\overline{x_{s}}]\ni a\mapsto  \big (\tau_\rT^{-1}\rP^{*}_{\rT}(t)+\rQ \big )g_{s\rT}(a),
\label{eq:maximumHamiltonian}
\end{equation}
for almost any $t$.
\begin{rem}\rm
We recall that in the Pontryagin Principle, the optimization is governed by the direct dependence on the control values $a$ (see \cite{FR}). $\fin$
\end{rem}
\begin{rem}\rm
The solution of the  functional minimization problem $\rW(\bV_{0})\doteq\min_{x_{s}\in \cU}\J _{\bV_{0}}(x_{s})$ (see \eqref{eq:minimization1})  is determined by \eqref{eq:maximumHamiltonian}.$\fin$
\end{rem}
\par
Then, from \eqref{eq:maximumHamiltonian} we have obtained an optimal control given by
\begin{equation}
x_{s}^{*}(t)=
\left \{
\begin{array}{ll}
\underline{x_{s}}, & \hbox{if $\rP^{*}_{\rT}(t)<-\tau_{\rT}\rQ$},\\ [.1cm]
\overline{x_{s}}, & \hbox{if $\rP^{*}_{\rT}(t)> -\tau_{\rT}\rQ$}. 
\end{array}
\right .
\label{eq:optimalcontrol}
\end{equation}
The time for which $\rP^{*}_{\rT}(t^{*})=-\tau_{\rT}\rQ$ is called a switching time.
\begin{rem}\rm We note that if $ \rP^{*}_{\rT}(t^{*})=-\tau_{\rT}\rQ $ holds for some 
$t^{*}$, the function
$$
a\mapsto \rH\big (\bV^{*}(t^{*}),\bP^{*}(t^{*}),a\big )
$$
is constant (see \eqref{eq:maximumHamiltonian}). Then, the value $x_{s}^{*}(t^{*})\in\big [\underline{x_{s}},\overline{x_{s}}\big ]$ is irrelevant.
\label{obser:irrelevancy}
\end{rem}
\begin{rem}\rm
Since $\rP^{*}_{\rT}(t_{f})=-\rS'\big (\rV_{\rT}^{*}(t_{f})\big )\in [-\rS_{+}',-\rS_{-}']$, if $\rS'_{+}<\tau _{\rT}\rQ$, we deduce that there exists  $\varepsilon >0$ such that 
\begin{equation}
x_{s}^{*}(t)\equiv \overline{x_{s}} \quad  \hbox{for $t_{f}-\varepsilon <t\le t_{f}$}.
\label{eq:finaloptimality-}
\end{equation}
Thus, near the final time $t_{f}$ the control value does not change.$\fin$
\end{rem}
\begin{theo}[Uniqueness of the optimal control] Assume 
\begin{equation}
\rS'_{+}<\tau_{\rT}\rQ
\label{eq:finalnochange}
\end{equation}
and the compatibility condition 
\begin{equation}
0<\ln \left (\dfrac{\rK-\rS'_{-}}{\rK-\tau_{\rT}\rQ}\right )^{\tau_{\rT}}<t_{f}, 
\label{eq:compatibility}
\end{equation}
where $\rK$ is a constant large enough such that
\begin{equation}
\rK>\max \left \{-\min_{t \in[0,\tau_{f}]} \rG_{\rE}'\big (\rV_{\rE}^{*}(t)\big )\rP_{\rI}^{*}(t)+\rV_{\rT}^{*}(t),\rS_{-}',\tau_{\rT}\rQ,\dfrac{\tau_{\rT}\rQ e^{\frac{t_{f}}{\tau_{\rT}}}-\rS'_{-}}{e^{\frac{t_{f}}{\tau_{\rT}}}-1}\right \}.
\label{eq:Klarge}
\end{equation}
Then the optimal control is the bang-bang type control  
\begin{equation}
x_{s}^{*}(t)=
\left \{
\begin{array}{ll}
	\underline{x_{s}}, & 0<t< t _{*},\\ [.1cm]
	\overline{x_{s}}, &   t _{*}<t< t_{f}.
\end{array}
\right .
\label{eq:uniqueoptimalcontrolcontrol}
\end{equation}
Moreover, the switching time $t ^{*}$ verifies the lower estimate  
\begin{equation}
0<t_{f}	-\ln \left (\dfrac{\rK-\rS'_{-}}{\rK-\tau_{\rT}\rQ}\right )^{\tau_{\rT}}<t _{*}. 
\label{eq:loweresimatetime}
\end{equation}
\label{theo:PTdecre}
\end{theo}
\proof The final problem on $\rP_{\rT}^{*}$ is
$$
\tau_{\rT}\derp{\rP_{\rT}^{*}}(t)=\rP_{\rT}^{*}(t)-\rG_{\rE}'\big (\rV_{\rE}^{*}(t)\big )\rP_{\rI}^{*}(t)-\rV^{*}_{\rE}(t), 
\quad \rP_{\rT}^{*}(t_{f})=-\rS'\big (\rV_{\rT}^{*}(t_{f})\big ).
$$
Then we get the inequality 
$$
\tau_{\rT}\derp{\rP_{\rT}^{*}}(t)> \rP_{\rT}^{*}(t)+\rK,\quad \hbox{ with }
\rK>-\max_{t \in[0,t_{f}]}\left (\rG_{\rE}'\big (\rV_{\rE}^{*}(t)\big )\rP_{\rI}^{*}(t)+\rV_{\rT}^{*}(t )\right ). 
$$
Therefore
$$
\left (\derp{\rP_{\rT}^{*}}(t)-\dfrac{\rP^{*}_{\rT}(t)}{\tau_{\rT}}\right )e^{-\frac{t}{\tau_{\rT}}}> \dfrac{\rK}{\tau_{\rT}}e^{-\frac{t}{\tau_{\rT}}}
\quad
\Leftrightarrow \quad \dfrac{d}{dt}\left (\big (\rP_{\rT}^{*}(t)+\rK\big )e^{-\frac{t}{\tau_{\rT}}}\right )>0,
$$
implies 
\begin{equation}
\big (\rP_{\rT}^{*}(t _{2})+\rK\big )e^{-\frac{t_{2}}{\tau_{\rT}}}> \big (\rP_{\rT}^{*}(t _{1})+\rK\big )e^{-\frac{t_{1}}{\tau_{\rT}}},
\label{eq:generalizedmononicity}
\end{equation}
whenever $t _{2}>t_{1}$. 
Next, by means of \eqref{eq:Klarge}, we prove that assumption \eqref{eq:finalnochange} implies the existence of a unique switching time $t_{*}$. Indeed, otherwise
$$
\rP_{\rT}^{*}(t)>-\tau_{\rE}\rQ\quad \hbox{for all $t $ in the interval $[0,t_{f}]$}
$$
(see \eqref{eq:finalnochange}). Then, from \eqref{eq:generalizedmononicity} we obtain
$$
\big (\rK-\rS_{-}'\big )e^{-\frac{t_{f}}{\tau_{\rT}}}>
\big (\rP_{\rT}^{*}(t_{f})+\rK\big )e^{-\frac{t_{f}}{\tau_{\rT}}}> \big (\rP_{\rT}^{*}(t)+\rK\big )e^{-\frac{t}{\tau_{\rT}}}> \big (\rK-\tau_{\rT}\rQ\big )e^{-\frac{t}{\tau_{\rT}}},
$$
whence by choosing $\rK>\tau_{\rT}\rQ$ we deduce $\rK>\rS_{-}'$ and 
$$
\big (\rK-\rS_{-}'\big )e^{-\frac{t_{f}-\sigma }{\tau_{\rT}}}>\rK-\tau_{\rT}\rQ\quad \hbox{for all  $t$ in the interval $[0,t_{f}]$}.
$$
Then, since  $\rS_{-}'\le \rS_{+}'<\tau_{\rT}\rQ$ (see \eqref{eq:finalnochange}) we may derive a contradiction by choosing $\widehat{t}<t_{f}$ given by 
$$
\big (\rK-\rS_{-}'\big )e^{-\frac{t_{f}-\widehat{t}}{\tau_{\rT}}}=\rK-\tau_{\rT}\rQ\quad
\Leftrightarrow \quad 
0<t_{f}-\ln \left (\dfrac{\rK-\rS'_{-}}{\rK-\tau_{\rT}\rQ }\right )^{\tau_{\rT}}=\widehat{t}<t_{f},
$$
provided
$$
0<\ln \left (\dfrac{\rK-\rS'_{-}}{\rK-\tau_{\rT}\rQ }\right )^{\tau_{\rT}}<t_{f} \quad \Leftrightarrow \quad \rK>\dfrac{\tau_{\rT}\rQ e^{\frac{t_{f}}{\tau_{\rT}}}-\rS'_{-}}{e^{\frac{t_{f}}{\tau_{\rT}}}-1},
$$
where we are assuming \eqref{eq:Klarge}.
\par
\noindent
Consequently, there exists a time $t_{*}\in [0,t_{f}[$, with $\rP^{*}_{\rT}(t_{*})=-\tau_{\rT}\rQ$.  In fact, in this case, we also may deduce from the contradiction
\begin{equation}
\big (\rK-\tau_{\rT}\rQ\big )e^{-\frac{\sigma _{2}}{\tau_{\rT}}}> \big (\rK-\tau_{\rT}\rQ\big )e^{-\frac{\sigma _{1}}{\tau_{\rT}}}>\big (\rK-\tau_{\rT}\rQ\big )e^{-\frac{\sigma _{2}}{\tau_{\rT}}}
\label{eq:contra}
\end{equation}
(see \eqref{eq:generalizedmononicity}) that 
there exists a unique time of change $t ^{*}<t_{f}$ with $\rP^{*}_{\rT}(t_{*})=-\tau_{\rT}\rQ$. 
\par
The above reasoning enables us to obtain a lower estimate of $t_{*}$. Indeed, by considering  \eqref{eq:generalizedmononicity}, we may construct the inequality 
$$
\big (\rK-\rS_{-}'\big )e^{-\frac{t_{f}-t_{*}}{\tau_{\rT}}}\ge
\big (\rP^{*}_{\rT}(t_{f})+\rK\big )e^{-\frac{t_{f}-t_{*}}{\tau_{\rT}}}\ge \rK-\tau_{\rT}\rQ,
$$
whence 
$$
0<t_{f}-\ln \left (\dfrac{\rK-\rS'_{-}}{\rK-\tau_{\rT}\rQ}\right )^{\tau_{\rT}}=\widehat{t}<t_{*}.
$$
\fineq
\begin{rem}\rm By means of the arguments of \eqref{eq:generalizedmononicity} and \eqref{eq:contra}, one proves
$$
\rP_{\rT}^{*}(t _{2})>\rP_{\rT}^{*}(t _{1}),\quad t_{f}\ge t _{2}>t _{1}\ge t_{*}.
$$
Then $\rP_{\rT}^{*}(t) $ is an increasing function on $[t_{*},t_{f}]$. We note that the upper bound 
$$
\rP_{\rT}^{*}(t)<-\tau_{\rT}\rQ \quad \hbox{if $0\le t <t _{*}$}
$$
holds.$\fin$
\end{rem}
\begin{rem}\rm The optimal control given in \eqref{eq:uniqueoptimalcontrolcontrol} admits other alternative representations. For instance, we have
\begin{equation}
x_{s}^{*}(t)=\dfrac{\overline{x_{s}}-\underline{x_{s}}}{2}\big (1+\hbox{\bf sign}(t-t_{*})\big )+\underline{x_{s}},\quad 0<t<t_{f}.
\label{eq:uniqueotimalcontrolcontrol2}
\end{equation}
\end{rem}
\par
\medskip
 To estimate  $t_{*}$ from the nonlinear equation
$$
\rP^{*}_{\rT}(t_{*})=-\tau_{\rT}\rQ
$$
is a very tedious task. An alternative argument is provided in the following result.
\begin{theo}[Characterization of the optimal switching time] Under the assumptions of Theorem \ref{theo:PTdecre}, the optimal time of change $t_{*}$ of the optimal control {\rm (}see \eqref{eq:optimalcontrol}{\rm )} is characterized by the condition
\begin{equation}
	\dfrac{1}{2}\dfrac{\partial }{\partial t_{*}}\int^{t_{f}}_{0}|\rV_{\rT}^{*}\big (\sigma ;t_{*}\big )|^{2}d \sigma-\rQ\big (g_{s\rT}\big (\overline{x_{s}}\big )-g_{s\rT}\big (\underline{x_{s}}\big )\big )+\rS'\big (\rV_{\rT}^{*}(t_{f};t_{*})\big)\dfrac{\partial \rV_{\rT}^{*}(t_{f};t_{*})\big)}{\partial t_{*}}=0,
	\label{eq:characterizationoptimaltimechange}	
\end{equation}
where we have denoting by $\rV_{\rT}^{*}(t;t_{*})$ to the optimal state .
\label{teo:caracterizacioninstanteotimocambio}
\end{theo}
\par
\noindent 
{\sc Proof} For every time $\widehat{t}_{*}\in [0,t_{f}]$ we consider a general bang-bang control 
\begin{equation}
x_{s}^{*}(t;\widehat{t}_{*})=
\left \{
\begin{array}{ll}
	\overline{x_{s}}, & 0\le t<\widehat{t}_{*},\\ [.1cm]
	\underline{x_{s}}, & \widehat{t}_{*}<t\le t_{f}.
\end{array}
\right .
\label{eq:bangbangcasioptimo}
\end{equation}
Then, by the above Theorem we have that $\widehat{t}_{*}=t_{*}\in\big (0,t_{f}\big )$ for the {\em optimal control} \eqref{eq:optimalcontrol}.
Let us denote by $\rVh_{\rT}^{*}(t;\widehat{t}_{*})$ the associate state corresponding to the control 	$x_{s}^{*}(t;\widehat{t}_{*})$. We define the function  
$$
\Phi (\widehat{t}_{*})=	\J (x_{s}^{*}(\cdot;\widehat{t}_{*})).
$$
Therefore, one has
$$
\min_{x_{s}\in \cU}\J (x_{s})=\J (x_{s}^{*})=\Phi (t_{*})=\min_{\widehat{t}_{*}\in [0,t_{f}]}\Phi (\widehat{t}_{*}).
$$
Thus, the optimal switching time $t_{*}$ must satisfies
$$
\left .\dfrac{d \Phi (\widehat{t}_{*})}{d \widehat{t}_{*}}\right |_{\widehat{t}_{*}=t_{*}}=0.
$$
$\fin$
\par
\medskip
Coming back to the equation of the $\rV_{\rE}$  {\em optimal potential}  
$$
\tau_\rE \derp{\rV_{\rE}^{*}}=-\rV_\rE^{*}+\rV_{\rE 0} + g_{s\rE}(\underline{x_{s}}),\quad 0\le t<t_{*}. 
$$
Straightforward computations yield
$$
\rV_{\rE}^{*}(t)=g_{s\rE}\big (\underline{x_{s}}\big )\left (1-e^{-\frac{t}{\tau_{\rE}}}\right )+\rV_{\rE 0},
\quad 0\le t<t_{*}.
$$
Then 
\begin{equation}
\begin{array}{ll}
\rV_{\rE}^{*}(t)&\hspace*{-.2cm}=g_{s\rE}\big (\underline{x_{s}}\big )\left (1-e^{-\frac{t}{\tau_{\rE}}}\right )+\rV_{\rE 0},\quad 0\le t\le t_{*},\\ [.2cm]
& \quad \left (\rV_{\rE 0}+g\big (\overline{x_{s}}\big )\right )\left (1-e^{-\frac{t-t_{*}}{\tau_{\rE}}}\right )+\left (
g\big (\underline{x_{s}}\big )\left (1-e^{-\frac{t_{*}}{\tau_{\rE}}}\right )+\rV_{\rE 0}\right )
e^{-\frac{t-t_{*}}{\tau_{\rE}}},\quad t_{*}<t\le t_{f}.
\end{array}
\end{equation}
$\fin$
\begin{center}

{\bf \sffamily Acknowledgements}
\end{center}
The research of the authors was partially supported by the project  PID2020-112517GB-I00 of the AEI (Spain). JID is also supported by the project MCIU/AEI/10.13039/-501100011033/FEDER, EU.

\par
\noindent
\vskip1cm
\flushleft{\small
\begin{tabular}{ll}\small
G. D\'{\i}az & J.I. Díaz \\
& Instituto Matemático Interdisciplinar (IMI)\\
Dpto. Análisis Matemático y Matem\'atica Aplicada & Dpto. Análisis Matemático y Matem\'atica Aplicada\\
Universidad Complutense de Madrid & Universidad  Complutense de Madrid\\
Parque de las Ciencias & Parque de las Ciencias \\
28040-Madrid. Spain & 28040-Madrid. Spain\\
{\tt gdiaz@ucm.es} & {\tt jidiaz@ucm.es}
\end{tabular}
 }

\begin{thebibliography}{99}

\bibitem[Am]{Am} Amann, H. (1990). \textit{Ordinary Differential Equations: An Introduction to Nonlinear Analysis}. Walter de Gruyter, Berlin.

\bibitem[An]{An} An der Heiden, U. (1980). \textit{Analysis of Neural Networks}. Lecture Notes in Biomathematics, \textbf{35}. Springer-Verlag, Berlin.

\bibitem[BP]{BP} Bressan, A. and Piccoli, B. (2007). \textit{Introduction to the Mathematical Theory of Control}. American Institute of Mathematical Sciences, Springfield, USA.

\bibitem[BCS]{BCS} Britton, N. F., Chaplain, M. A. J., and Skevington, S. M. (1996). The role of N-methyl-D-aspartate (NMDA) receptors in wind-up: A mathematical model. \textit{J. Math. Appl. Med. Biol.}, \textbf{13}(2), 193-205.

\bibitem[BS]{BS} Britton, N. F., and Skevington, S. M. (1988). A mathematical model of the Gate Control Theory of Pain. \textit{J. Bio. Syst.}, \textbf{3}(4), 1119-1124.

\bibitem[CD]{CD} Casal, A. C., and Díaz, J. I. (2005). On the principle of pseudo-linearized stability: Applications to some delayed nonlinear parabolic equations. \textit{Nonlinear Anal.}, \textbf{63}(8), 997-1007.

\bibitem[FR]{FR} Fleming, W., and Rishel, R. (1975). \textit{Deterministic and Stochastic Optimal Control}. Springer, New York.

\bibitem[MW1]{MW1} Melzack, R., and Wall, P. D. (1965). Pain mechanisms: A new theory. \textit{Science}, \textbf{150}(3699), 971-979.

\bibitem[MW2]{MW2} Melzack, R., and Wall, P. D. (1982). \textit{The Challenge of Pain}. Penguin, Harmondsworth, UK.

\bibitem[Na]{Na} Nathan, P. W. (1976). The gate control theory of pain: A critical review. \textit{Brain}, \textbf{99}(1), 123-158.

\bibitem[NR]{NR} Nathan, P. W., and Rudge, P. (1974). Testing the gate control theory of pain in man. \textit{J. Neurol. Neurosurg. Psych.}, \textbf{37}(12), 1366-1372.

\bibitem[Tr]{Tr} Trélat, E. (2008). \textit{Contrôle Optimal: Théory et Applications}. Vuibert, Collection Mathématiques Concrètes, Paris.

\bibitem[WC]{WC} Wilson, H. R., and Cowan, J. D. (1972). Excitatory and inhibitory interactions in localized populations of model neurons. \textit{Biophys. J.}, \textbf{12}(1), 1-24.

\end{thebibliography}
\end{document}